\documentclass[pdflatex,sn-mathphys-num]{sn-jnl}% Math and Physical Sciences Numbered Reference Style
\usepackage{graphicx}%
\usepackage{multirow}%
\usepackage{amsmath,amssymb,amsfonts}%
\usepackage{amsthm}%
\usepackage{mathrsfs}%
\usepackage[title]{appendix}%
\usepackage{xcolor}%
\usepackage{textcomp}%
\usepackage{manyfoot}%
\usepackage{booktabs}%
\usepackage{algorithm}%
\usepackage{algorithmicx}%
\usepackage{algpseudocode}%
\usepackage{listings}%
\renewcommand{\i}{\mathbf{i}\,}
\DeclareMathOperator{\Tr}{Tr}

\theoremstyle{thmstyleone}%
\theoremstyle{thmstyletwo}%

\theoremstyle{thmstylethree}%

\begin{document}

\title[Spectra of random graphs with discrete scale invariance]{Spectra of random graphs with discrete scale invariance}


\author*[1,2,3]{\fnm{Alessio} \sur{Catanzaro}}\email{alessio.catanzaro@imtlucca.it}\equalcont{These authors contributed equally to this work.}

\author[4]{\fnm{Rajat} \sur{Subhra Hazra}}\email{r.s.hazra@math.leidenuniv.nl}
\equalcont{These authors contributed equally to this work.}

\author[1,2]{\fnm{Diego} \sur{Garlaschelli}}\email{diego.garlaschelli@imtlucca.it}
\equalcont{These authors contributed equally to this work.}

\affil*[1]{\orgdiv{Network Unit}, \orgname{IMT Lucca}, \orgaddress{ \city{Lucca},  \country{Italy}}}

\affil[2]{\orgdiv{Lorentz Institute for Theoretical Physics}, \orgname{Leiden University}, \orgaddress{ \city{Leiden}, \state{ZH}, \country{The Netherlands}}}

\affil[3]{\orgdiv{Dipartimento di Fisica}, \orgname{Università di Palermo}, \orgaddress{\city{Palermo},  \country{Italy}}}

\affil[4]{\orgdiv{Mathematical Institute}, \orgname{Leiden University}, \orgaddress{ \city{Leiden}, \state{ZH}, \country{The Netherlands}}}


\abstract{Random graphs defined by an occurrence probability invariant under node aggregation have been identified recently in the context of network renormalization. The invariance property requires to draw edges with a specific probability that, in the annealed case, depends on a necessarily infinite-mean node fitness. The diverging mean determines many properties that are uncommon in models with independent edges, but at the same time widespread in real-world networks. Here we focus on leading eigenvalues and eigenvectors of the adjacency matrix of the model, where the $n$ nodes are assigned a Pareto($\alpha$)-distributed fitness with $0<\alpha<1$. We find that the leading eigenvalues are of order $\sqrt{n}$, alternate in sign and are located at the intersection between the real axis and a logarithmic spiral in the complex plane, which we characterize analytically in terms of the Gamma function. We also calculate the associated eigenvectors, finding complex-valued scaling exponents and log-periodicity, which are signatures of discrete scale-invariance. In contrast with typical finite-rank behaviour of random graphs, we find that a growing number of the leading eigenvalues emerges from the bulk, whose edge extends up to order $\sqrt{n}$ and therefore reaches the same scale as that of the structural eigenvalues.}




\maketitle

\section*{Introduction}

Many models have been proposed to reproduce the heterogeneous structure ubiquitously observed in real-world networks. Among these, a wide class of random graph models assigns a weight (or fitness, or hidden variable) $x_i$ to each vertex $i$ ($i=1,n$ where $n$ is the number of nodes) from a usually fat-tailed distribution $\rho(x)$, and connects pairs of vertices $i$ and $j$ with a probability $p_{ij}$ that is a certain positive function $f(x_i,x_j)$~\cite{caldarelli2002scale,garlaschelli2004fitness,van2017random,van2024random}.

Recently, in the context of network renormalization~\cite{gabrielli2025network, catanzaro2026renormalization}, a specific such model that complements the heterogeneity of nodes with an invariance principle under graph coarse-graining  has been identified~\cite{garuccio2023multiscale,lalli2024geometry, avena2022inhomogeneous, milocco2024multi, catanzaro2026clusteringgeometrysparsenetworks}. 
In this so-called (annealed) Multi-Scale Model (MSM), one may group nodes into equally sized `supernodes', each inheriting the sum of the fitness of its members, and impose that both the connection probability $p_{ij}=f(x_i,x_j)$ and the fitness distribution $\rho(x)$ remain unchanged in their functional form, while their parameters follow exact renormalization rules. 
In other words, aggregating nodes into supernodes according to any homogeneous partition -- and connecting two supernodes whenever there is a connection between any of their constituent nodes -- defines a coarse-grained version of the graph, governed by an occurrence probability that has the same functional form as that of the original graph, up to a rescaling of the parameters~\cite{garuccio2023multiscale,avena2022inhomogeneous}.

This renormalization scheme does not rely on any underlying geometry or distance notion between nodes. 
This stands in contrast to geometric renormalization approaches \cite{garcia2018multiscale, boguna2021network}, where the nodes to be merged need be spatially proximate, and a preliminary definition of metric distance is needed.
Rather, when nodes are merged into supernodes of equal cardinality, the invariance requirement in the MSM is closer in spirit to the concept of discrete scale invariance~\cite{logperiodic1}, whose consequences for the graph properties have not been investigated so far and are one of our main focuses here.

In particular, we are interested in the spectral properties~\cite{van2023graph} of the MSM, which turn out to be quite unique due to the peculiar structure of the fitness. Indeed, the requirement that the node fitness is a positive random (`annealed') variable that is stable under addition implies that it has an $\alpha$-stable distribution with infinite mean ($0<\alpha<1$)~\cite{garuccio2023multiscale,avena2022inhomogeneous}, where classical random matrix results~\cite{benaych2016lectures, bai2010spectral, van2023graph, farkas2001spectra, poley2024eigenvalue, baron2022eigenvalues, chung2003spectra} break down. 
When a symmetric random matrix is decomposed into a low-rank informative part plus noise, strong signals can produce eigenvalue outliers outside the bulk spectrum, as first shown in \cite{baik2005phase}. In networks, such outliers reveal structures like communities or hubs \cite{nadakuditi2013spectra, van2023graph}. %However the heavy-tailed fitness distribution poses a challenge: hubs generate multiple large eigenvalues, but standard random matrix theory fails due to diverging moments. 
While infinite-mean regimes have been studied for L\'evy matrices \cite{cizeau1994theory,benaych2014central,arous2008spectrum}, random graphs add further randomness from adjacency realizations, not captured by those approaches.
Our goal is to shed light on the surprising behaviour of the eigenmodes of the adjacency matrix in this infinite-mean regime, which follows from discrete node aggregation invariance. From a network-science perspective, the leading eigenvalues and eigenvectors are useful because they describe the large-scale structure created by the node fitnesses. In an actual graph realization this structure is mixed with the randomness of the edge variables. The leading spectral modes help identify which part of the spectrum reflects the multi-scale organization of the model, and which part belongs to the random bulk. In standard finite-mean hidden-variable models, only finitely many such structural modes are usually expected to stand out. In the MSM, by contrast, the number of visible modes grows slowly with the system size. This makes the leading spectrum a possible signature of discrete scale invariance and of the log-periodic behaviour generated by repeated coarse-graining. More generally, our analysis suggests a way to study random graph spectra in models where the expected adjacency matrix is not finite-rank, but still has only a slowly growing number of relevant modes.

\section*{Results}
\subsection*{Model definition and matrix decomposition}
The annealed MSM is an inhomogeneous random graph on $n$ nodes, each node $i$ being assigned a positive weight $x_i$, sampled independently from an $\alpha$-stable distribution $\rho_\alpha(x)$, hence with power-law tails decaying as $x^{-1-\alpha}$ for $0<\alpha<1$~\cite{garuccio2023multiscale}. 
In this regime of $\alpha$, every moment of $\rho_\alpha(x)$ diverges.
Given $x_i$ and $x_j$, the probability that an edge exists between $i$ and $j$ is \begin{equation}\label{def:f_n}
    p_{ij}= f_n(x_i, x_j)=1-\exp(-\varepsilon_n x_i x_j),\quad i,j=1,n
    \end{equation} 
if $i\ne j$ and $p_{ii}=0$ if $i=j$, where the parameter $\varepsilon_n$ may be chosen to scale with $n$ in a desired way. 
The exponential form in Eq.~\eqref{def:f_n} is used due to its compatibility with the MSM coarse-graining rule. Indeed, if $I$ and $J$ are two supernodes with aggregated fitnesses
\[
X_I=\sum_{i\in I}x_i,\qquad X_J=\sum_{j\in J}x_j 
\]
(where $i \in I$ denotes that node $i$ belongs to the supernode $I$), and if two supernodes are connected whenever at least one microscopic edge between them is present, then
\[
1-\prod_{i\in I,j\in J}(1-p_{ij})
=
1-\exp\left(-\varepsilon_n\sum_{i\in I,j\in J}x_i x_j\right)
=
1-\exp(-\varepsilon_n X_I X_J).
\]
Thus the functional form of the connection probability is preserved exactly under homogeneous aggregation.
Other choices, such as $\min\{\varepsilon_n x_i x_j,1\}$ or $\varepsilon_n x_i x_j/(1+\varepsilon_n x_i x_j)$, have the same first-order behaviour when $\varepsilon_n x_i x_j$ is small, but they do not satisfy the above exact no-edge multiplicative identity. We therefore expect some qualitative features, such as the $\sqrt n$-scale of the leading modes and the appearance of power-law eigenfunctions, to be robust for a broader class of saturating link functions. However, the exact Gamma-function constants and the explicit logarithmic spiral are specific to the exponential connection probability, because they follow from the Laplace transform calculation used below.

Under an arbitrary homogeneous partition whereby the $n$ nodes are aggregated into supernodes (each with the same number of participating nodes), the form of both $p_{ij}=f_n(x_i,x_j)$ and $\rho_\alpha(x)$ is unchanged, up to a parameter shift in the latter -- reflecting the fact that the fitness of a supernode is the sum of the fitnesses of its constituent nodes~\cite{garuccio2023multiscale,avena2022inhomogeneous}.
Since $\rho(x)$ for positive $\alpha$-stable variables is not known in closed-form (except for $\alpha=1/2$), for analytical tractability we follow~\cite{avena2022inhomogeneous} and consider a pure Pareto($\alpha$) distribution with the same tail as an $\alpha$-stable distribution:
\begin{equation}
    \rho_\alpha(x) = \alpha~ x^{-1-\alpha}\quad 0<\alpha<1
\end{equation} for $x\ge 1$, and $\rho_\alpha(x) =0$ otherwise.
We also consider $\varepsilon_n = n^{-1/\alpha}$ as a scaling that allows the graph to be arbitrarily large while remaining sparse, allowing precise asymptotic calculations for large $n$~\cite{avena2022inhomogeneous}. In this regime, the tail of the expected degree distribution $P(k)$ decays as $k^{-2}$ irrespective of the value of $\alpha$, the expected link density is of order $\ln n/n$, while the average local clustering coefficient remains finite~\cite{garuccio2023multiscale,avena2022inhomogeneous} -- which are all widespread properties found in real-world networks.

\subsection*{Realized versus expected matrix} 
We denote the realized adjacency matrix of the graph as $\mathbf{A}$, where $A_{ij}$ is a Bernoulli variable: $A_{ij}=1$ with probability $p_{ij}$ and $0$ otherwise. Note that $\mathbf{A}$ is symmetric and $A_{ii}=0$. By construction, the expectation of $\mathbf{A}$ (conditioning on the realized weights) is the matrix $\mathbf{P} = \mathbb{E}[\mathbf{A} | \{x_i\}]$ with entries $P_{ij} = p_{ij} = f_n(x_i,x_j)$ and captures the deterministic, structural signal imposed by the weights, while the actual $\mathbf{A}$ can be seen as a noisy realization of that signal.

We use the term ``annealed'' for the sampling of the fitness variables. Conditional on the realized weights $\{x_i\}$, however, the edge variables are independent Bernoulli variables and the expected matrix is $\mathbf P=\mathbb E[\mathbf A\mid \{x_i\}]$. Thus the spectral comparison between $\mathbf A$, $\mathbf P$, and $\mathbf H=\mathbf A-\mathbf P$ is naturally quenched with respect to the edge randomness. If the random weights are replaced by a deterministic sequence with the same regularly varying order statistics, for instance \(x_i=(n/i)^{1/\alpha}\) after ordering, the leading asymptotic predictions remain the same, while sample-to-sample fluctuations of the outlier locations are suppressed. This is also illustrated in Fig.~\ref{fig:lambda_max}, where the deterministic and sampled Pareto weights display the same $\sqrt n$ scaling.

\begin{figure}[t]
     \centering
     \includegraphics[width=\linewidth]{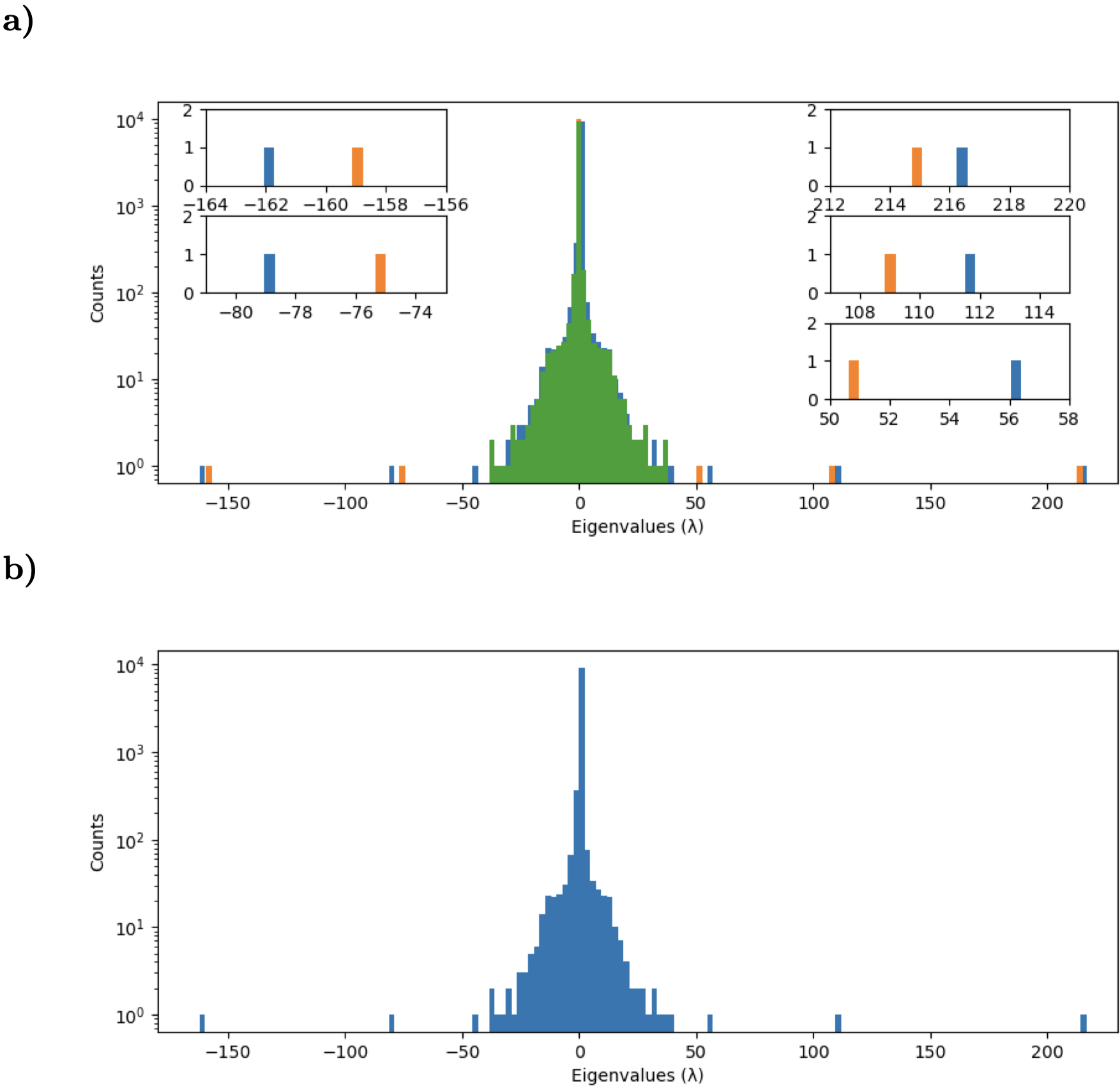}
     \caption{\textbf{Full eigenvalue spectrum of the adjacency matrix $\mathbf{A}$ and its decomposition into outlier and bulk components for $n=10^4$ and $\alpha = 0.5$.}
\textbf{a)} Histogram (blue) of the eigenvalues of $\mathbf{A}$. We see separation into atomic delta densities and a central, continuous bulk.
\textbf{b)} Separation of the spectrum into two regimes, the outliers being those of $\mathbf{A}$ in blue, and of $\mathbf{P}$ in orange. We also superimpose in green the eigenvalue histogram of $\mathbf{H}$, that replicates the bulk behaviour of $\mathbf{A}$.}
     \label{fig:fullspectrum}
 \end{figure}

Indeed, a central idea in our analysis is to decompose the adjacency matrix into the sum of an expected component and a fluctuation component, and then use Random Matrix Theory (RMT) intuition to separate the spectrum into bulk and outlier components. Let us write: \begin{equation}\label{eq:A-P-H}
\mathbf{A} = \langle \mathbf{A} \rangle + (\mathbf{A} - \langle \mathbf{A} \rangle) = \mathbf{P} + \mathbf{H}, 
\end{equation}  
where $\mathbf{H} = \mathbf{A} - \mathbf{P}$ is the zero-mean noise matrix. By construction, $\mathbf{P}$ is a dense matrix (encoding for the structure of the network) with a highly skewed pattern due to the heavy-tailed ${x_i}$, whereas $\mathbf{H}$ captures the randomness of edge realizations. 
If the matrix $\mathbf{P}$ were of small, finite rank and $\mathbf{H}$ were a `typical' random matrix with a bounded spectral distribution, then standard results show that each eigenvalue of $\mathbf{P}$ above a certain magnitude would detach from the continuous spectrum of $\mathbf{H}$. This phenomenon, rigorously characterized in spiked random matrix theory, underlies applications like spectral clustering \cite{zhang2014spectra, krzakala2013spectral, bordenave2015non}. In our case, $\mathbf{P}$ is not finite-rank: as we show in Supplementary Information (SI), we find evidence of an `effective' rank of order $\ln n$. 
Classical BBP-type results concern fixed finite-rank perturbations of random matrices. Our matrix $\mathbf P$ is not finite-rank, but its eigenvalues decay rapidly enough that only a slowly growing number of modes remains visible above the noise level. This places the MSM closer to recent work on growing-rank or mesoscopic perturbations, where the rank of the deterministic component may increase with $n$ while remaining $o(n)$~\cite{afanasiev2025asymptotic,huang2018mesoscopic}. We do not use these results directly, since our noise matrix has a heavy-tailed variance profile. The comparison is meant only to emphasize that the relevant perturbation is not finite-rank, but has an effective rank that grows sublinearly with $n$.
%In this sense, we are close in spirit to recent results such as Ref.~\cite{afanasiev2025asymptotic}, that are hinting at the fact that generalized BBP transitions could be happening even when the rank of the perturbation matrix is growing as $o(n)$. 
An illustration of the typical spectra of $\mathbf{A}$ and $\mathbf{P}$ is provided in Fig.~\ref{fig:fullspectrum}.
Our main goal is the characterization of these spectra and the associated eigenvectors, with a specific attention to the outlier modes carrying  structural information.

\subsubsection*{Integral eigenvalue equation} We are interested in the solution of the eigenmode equation~\cite{van2023graph}
\begin{equation}
\mathbf{P}v_k=\lambda_k v_k \quad k=1,n
\label{eq:eigen0}
\end{equation}
with $|\lambda_1|\ge |\lambda_2|\ge\dots\ge|\lambda_n|$. Denoting the $i$-th entry of each eigenvector $v_k$ as $v_k^{(i)}$, we can rewrite Eq.\eqref{eq:eigen0} as
\begin{equation}
\sum_{j=1}^n p_{ij}v^{(j)}_k=\lambda_k v^{(i)}_k \quad i=1,n\quad k=1,n  
\label{eq:eigen-}
\end{equation}
Note that, since $\mathbf{P}$ is real and symmetric, its eigenvalues are real.
Also note that any two vertices $i,j$ with the same realized weight $x_i=x_j$ are statistically equivalent, as their (expected) properties in the graph depend only on the weight $x$. Similarly, if $x_i=x_j$ then the $i$-th and $j$-th entries of each eigenvector $v_k$ will be equal: $v_k^{(i)}=v_k^{(j)}$. This allows us to write $v_k^{(i)}=\nu_k(x_i)$ where $\nu_k(x)$ must be a continuous function of the weight $x$.
For large enough $n$, Eq.\eqref{eq:eigen-} can be recast as the integral expression
%\begin{equation}
    %n\int_0^{+\infty}dy~ \rho_\alpha(y)~f_n(x,y)~\nu_k(y)=\lambda_k~\nu_k(x),
    %\label{eq:eigen1}
%\end{equation} %or equivalently 
\begin{equation}
    \mathbf{\Pi}_{\alpha,n}[\nu_k(x)]=\frac{\lambda_k}{n}~\nu_k(x).
    \label{eq:eigen2}
\end{equation} where $\nu_k(x)$ is an eigenfunction (with eigenvalue $\lambda_k/n$) of the integral operator
\begin{equation}
\mathbf{\Pi}_{\alpha,n}[g(x)]\equiv\int_0^{+\infty}dy~ \rho_\alpha(y)~(1-e^{-\varepsilon_n x y})~g(y),
    \label{eq:operator}
\end{equation} corresponding to the eigenvector $v_k$ (with eigenvalue $\lambda_k$) of the original matrix $\mathbf{P}$. A natural normalization for $\nu_k(x)$ will emerge from the calculations that follow.
If an eigenfunction $\nu_k(x)$ is found, a convenient way to go back to the actual eigenvector entry $v_k^{(j)}$ is the following. We first order all vertices so that $x_1\le x_2\le\dots\le x_n$.
Then, since a pure Pareto($\alpha$) variable can be represented as a uniform random variable,  raised to the power $-1/\alpha$, in the interval $[0,1]$, we arrive at the relation $x_j=(n/j)^{1/\alpha}$ (we will use this relation as an alternative, deterministic assignment of the fitness of nodes to produce a `noise-free' power law distribution even for finite $n$). Finally,
\begin{equation}
v_k^{(j)}=\nu_k(x_j)=\nu_k\left(\frac{n^{1/\alpha}}{j^{1/\alpha}}\right)\quad j=1,n\quad k=1,n.
\label{eq_mapping}
\end{equation}

\subsection*{Largest eigenvalue and principal eigenvector} 
\begin{figure}[b]
    \centering
    \includegraphics[width=\linewidth]{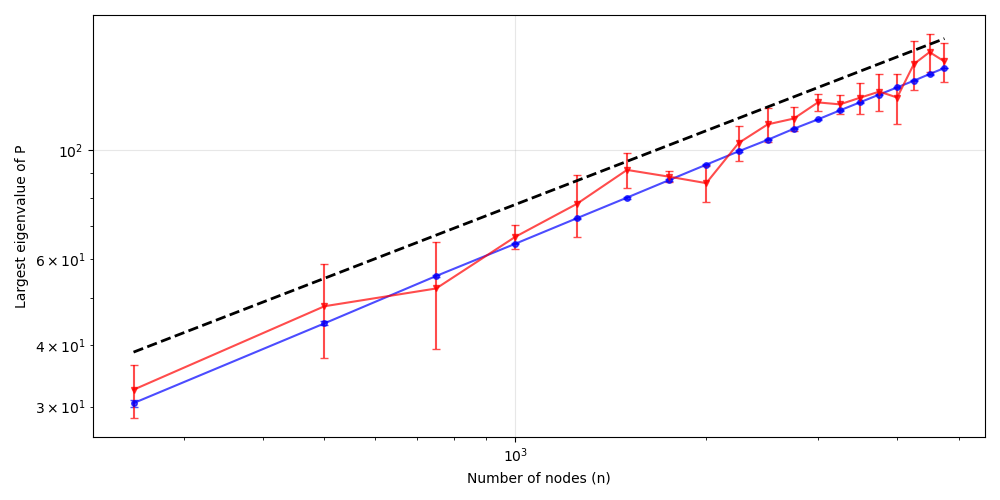}
    \caption{\textbf{Scaling of the largest eigenvalue $\lambda_1$ of $\mathbf{P}$ with the number of nodes $n$ for $\alpha=0.5$.}
The largest eigenvalue $\lambda_1$ is shown as a function of $n$ for tail index $\alpha$, using two assignments of the node weights $\{x_i\}_{i=1}^n$: weights sampled independently and identically distributed from a Pareto distribution (red circles, averaged over 5 independent samples) and weights assigned deterministically as $x_i=(n/i)^{1/\alpha}$ (blue circles). Both are compared with the analytic prediction of Eq.~\eqref{eq:eigenv1} (black solid line). The close agreement across both weight assignments confirms the $\sqrt{n}$ scaling.}
    \label{fig:lambda_max}
\end{figure}

We start with the principal component $k=1$. From Perron-Frobenius theorem, $\lambda_1$ is positive and $v_1$ has all entries of the same sign. 
We will set the normalization of $v_1$ in such a way that all its entries are positive, so that $\nu_1(x)$ is a positive-valued function. 
When evaluating the operator in Eq.~\eqref{eq:operator} on the eigenfunction, the multiplication of $\rho_\alpha(x)$ by $\nu_1(x)$ formally generates a modified, effective distribution of the weights. We then look for an eigenfunction such that $\nu_1(x)\rho_\alpha(x)=\rho_{\beta}(x)$ is again a Pareto pdf with a new index $\beta$. This is obtained for the eigenfunction \begin{equation}
    \nu_1(x)= \frac{\beta}{\alpha}~x^{\alpha-\beta},\quad \beta\in\mathbb{R^+}.
    \label{eq:nureal}
\end{equation} 
Setting $g(x)=\nu_1(x)$ into Eq.~\eqref{eq:operator}, we get \begin{equation}
    \mathbf{\Pi}_{\alpha,n}[\nu_1(x)] =1-\varphi_{\beta}(\varepsilon_n x)
      \label{eq:eigen3}
\end{equation} where we have introduced the Laplace Transform (LT) $\varphi_{\alpha}$ of the Pareto-$\alpha$ density, defined as
\begin{equation}
    \varphi_{\alpha}(t)\equiv\int_0^{+\infty}dx~ \rho_\alpha(x)~e^{-t x }=\alpha\int_1^{+\infty}dx~ x^{-1-\alpha}~e^{-t x },
\end{equation}
and we have exploited the normalization property 
\begin{equation}
    \int_0^{+\infty}dy~ \rho_\alpha(y)~\nu_1(y)=\int_0^{+\infty}dy~ \rho_{\beta}(y)=1.
    \label{eq:norm}
\end{equation}
The above expression sets our normalization for $\nu_1(x)$, which is the analogue of the $\ell_1$ norm $\sum_{i=1}^n |v_1^{(i)}|=1$ for the original eigenvector $v_1$.
An important asymptotic property of $\varphi_{\beta}(t)$ (see SI) is that, for $0\le\beta\le1$, \begin{equation}
1-\varphi_{\beta}(t)~\sim~ t^\beta ~\Gamma(1-\beta),\qquad t\to0^+.
\label{eq:LT}
\end{equation} This allows us to confirm that $\nu_1(x)$ is asymptotically an eigenfunction, for a suitable choice of $\beta$. Indeed, applying Eq.~\eqref{eq:LT} into Eq.~\eqref{eq:eigen3}, we can recast Eq.~\eqref{eq:eigen2} as
\begin{equation}
\mathbf{\Pi}_{\alpha,n}[\nu_1(x)]\sim (\varepsilon_nx)^\beta ~\Gamma(1-\beta)
    =\frac{\lambda_1}{n}~\frac{\beta}{\alpha}~x^{\alpha-\beta}
      \label{eq:eigen_k1}
\end{equation} (provided we identify $\beta\equiv\alpha-\beta$, i.e.  $\beta=\alpha/2$) to find
\begin{equation}
    \lambda_1=\frac{\alpha}{\beta}~\Gamma\left(1-\beta\right)n~\varepsilon_n^{\beta}=-\alpha~\Gamma\left(-\frac{\alpha}{2}\right)n^{1/2}.
    \label{eq:eigenv1}
\end{equation}
This square-root growth of $\lambda_1$  can be confirmed via actual computations
on the matrix $\mathbf{P}$ itself (see Fig.~\ref{fig:lambda_max}).
Using Eq.~\eqref{eq_mapping}, we also obtain all the eigenvector entries as
\begin{equation}
v_1^{(j)}=\nu_1(x_j)=\frac{1}{2}x_j^{\alpha/2}=\frac{1}{2}\left(\frac{n}{j}\right)^{1/2}\quad j=1,n.
\label{eq_principal}
\end{equation}
This is also confirmed numerically later on.

\subsection*{Next eigenvalues} 
\begin{figure}[h]
    \centering \includegraphics[width=\linewidth]{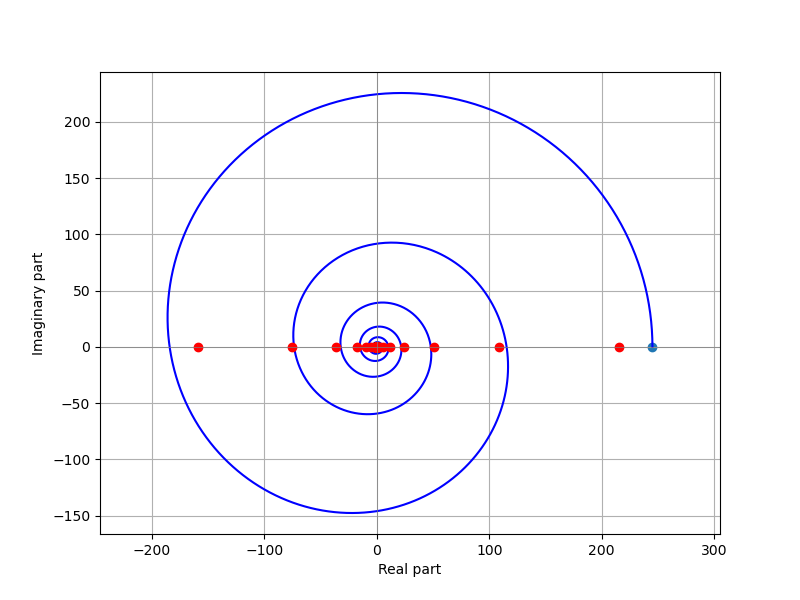} 
    \caption{\textbf{Eigenvalues of the expected adjacency matrix $\mathbf{P}$ lie at the intersections of the real axis with a logarithmic \emph{spira mirabilis} in the complex plane.}
The actual real-valued outlier eigenvalues of $\mathbf{P}$ (red dots) are shown together with the logarithmic spiral $\sigma^+(\omega)$ (blue curve) defined in Eq.~\eqref{eq:spiral}, plotted parametrically in the complex plane as a function of $\omega\ge0$. Parameters: $n=10^4$, $\alpha=0.5$.}
    \label{fig:spiral_05}
\end{figure}

We now turn $k>1$.
Let us preliminarily introduce auxiliary complex-valued test functions $\nu^\pm_k(x)$ of the form still given by Eq.~\eqref{eq:nureal}, where the role of the exponent $\beta$ is now taken by its complex counterparts $\beta^{\pm}_k \in \mathbb{C}$:   \begin{equation}
    \nu^\pm_k(x)= \frac{\beta_k^\pm}{\alpha}~x^{\alpha-\beta_k^\pm},\quad \beta_k^\pm\equiv\gamma_k\pm \i\omega_k\in\mathbb{C}
\end{equation}
where $\gamma_k=\Re[\beta_k^\pm]\in\mathbb{R}$ and $\pm \omega_k=\Im[\beta_k^\pm]\in\mathbb{R}$ are the real and imaginary parts of $\beta_k^\pm$, respectively. Note that $\beta_k^+$ and $\beta_k^-$ are complex conjugate.
Just like $\rho_{\beta}(x)\equiv \nu_1(x)\rho_\alpha(x)$ is a Pareto pdf with modified (real) index $\beta$, which has the same tail behaviour of a stable distribution with index $\beta$,  $\rho_{\beta_k^{\pm}}(x)\equiv \nu^\pm_k(x)\rho_\alpha(x)$ is formally a generalized Pareto distribution with the same `tail behaviour' as the recently defined stable distributions with complex-valued index $\beta_k^{\pm}$~\cite{alexeev2022stable,alexeev2023stable}. 
As we now show, these complex-valued scaling exponents generate log-periodicity, as typically found in systems with discrete scale invariance~\cite{logperiodic1}.
Indeed, in SI we derive the following expression for the eigenfunctions $\nu_{k}(x)$ for $k\ge1$:
\begin{equation}
\nu_k(x)=\frac{1}{2}\left[\Gamma(1-\beta_k^-)~\nu^+_k(x)+\Gamma(1-\beta_k^+)~\nu^-_k(x)\right]
%=(\varepsilon_nx)^{\gamma}(\varepsilon_nx)^{\pm i\omega} ~\Gamma(1-\gamma\mp i\omega_k^\pm)
      \label{eq:eigenf>}
\end{equation} 
with $\gamma_k=\alpha/2$ (for all $k$) and associated eigenvalues
\begin{equation}  \label{eq:eigenvk}  
     \lambda_k = (-1)^k~\alpha \left| {\Gamma\left(-\frac{\alpha}{2}+\i\omega_k\right)}\right|~ n^{1/2},
\end{equation}
where the values of $\omega_k$ 
are defined implicitly via
\begin{equation}
\arg{\left[{\Gamma\left(-\frac{\alpha}{2}+\i\omega_k\right)}\right]}=
    \frac{\omega_k}{\alpha}\ln n-k\pi,\quad k\in\mathbb{N}_0.
\label{eq:admissible}
\end{equation}
Note that $k=1$ and $\omega_1=0$ is a solution, coinciding with what we have already found in Eqs.~\eqref{eq:nureal} and~\eqref{eq:eigenv1} for the principal eigenmode. 
Graphically, $\lambda_1$ and all the other eigenvalues $\lambda_k$ ($k>1$) correspond to the intersections between the real axis and the two logarithmic spirals 
\begin{equation}
\sigma^\pm(\omega)
\equiv
-\alpha\,\Gamma\left(-\frac{\alpha}{2}\mp i\omega\right)
n^{1/2}
\exp\left(\pm \i\frac{\omega}{\alpha}\ln n\right).
\label{eq:spiral}
\end{equation}
defined parametrically in the complex plane for any real value $\omega\ge 0$ (see SI). 
In particular, plotting $\Im[\sigma^\pm(w)]$ versus $\Re[\sigma^\pm(w)]$ generates two spirals, both starting at $\lambda_1>0$ for $\omega=0$ and then shrinking clockwise ($\sigma^-$) and counterclockwise ($\sigma^+$) for $\omega>0$, crossing the real line at the next eigenvalues $\lambda_2$, $\lambda_3$, $\dots$ which alternate in sign.
The admissible values $\omega_k$ defined by Eq.~\eqref{eq:admissible} are precisely those realizing the intersections. 
A comparison between these intersections (for $\sigma^+$) and the actual eigenvalues computed directly on $\mathbf{P}$ for given  $n$ and $\alpha$ is illustrated in Fig.~\ref{fig:spiral_05} and in SI. 
Unfortunately, calculating the exact values $\omega_k$ defined by the non-invertible Eq.~\eqref{eq:admissible} is not possible for $k>1$. 
However, in SI we derive the following approximate expression (valid for $1<k\lesssim\ln n$):
\begin{equation}
    \omega_k \approx \alpha\frac{k\pi +\phi_\alpha}{\ln n},\quad \phi_\alpha\equiv\arg{\Gamma[-\alpha/2+\i\omega_\alpha]}
    \label{sollazzo}
\end{equation} 
where $\omega_\alpha\equiv\sqrt{\frac{\alpha}{2}(\frac{1}{\gamma}-\frac{\alpha}{2})}$ ($\gamma\approx 0.5772$ being the Euler-Mascheroni constant).
Equation~\eqref{sollazzo}  will prove useful to determine how fast $|\lambda_k|$ decays as $k$ grows. 

%Note that $k=1$ yields $\omega_1=0$, retrieving Eqs.~\eqref{eq:eigenv1}. In SI we show that the other values of $\omega_k$ satisfying \eqref{eq:admissible}, are described by the linear relation 

\subsection*{Next eigenvectors} 
\begin{figure}[h]
    \centering
    \includegraphics[width=\linewidth]{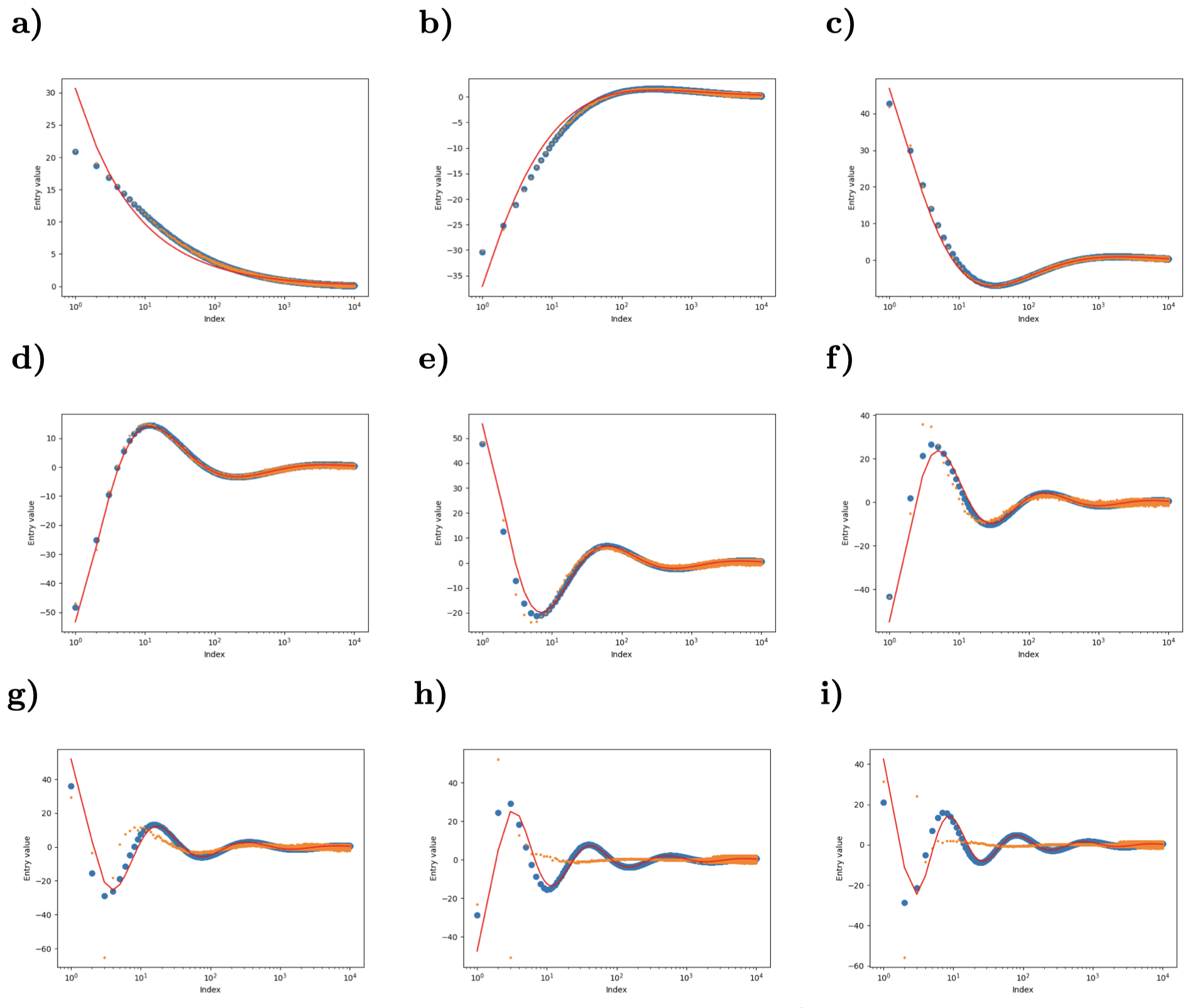}
    \caption{\textbf{Plot of eigenvector entries $v_k(j)$ versus $j$ for $n = 10^4$ and $\alpha = 0.5$.} \textbf{a)} First eigenvector. \textbf{b)} Second eigenvector. \textbf{c)} Third eigenvector. \textbf{d)} Fourth eigenvector. \textbf{e)} Fifth eigenvector. \textbf{f)} Sixth eigenvector. \textbf{g)} Seventh eigenvector. \textbf{h)} Eight eigenvector. \textbf{i)} Ninth eigenvector. Blue dots correspond to the actual eigenvectors entries of $\mathbf{P}$, while red solid lines represent the corresponding predicted entries from \eqref{eq:eigenvec}. Orange dots are the entries of the eigenvectors of $\mathbf{A}$, and we can see that in the outlier regime they are well predicted by $\mathbf{P}$ and hence the ones we obtained analytically. We also notice a breakpoint in $k$ after which the eigenvalues of $\mathbf{A}$ are no more described by those of $\mathbf{P}$, which are still predicted by \eqref{eq:eigenvec}. The agreement between the eigenvectors of $\mathbf{A}$ and those of $\mathbf{P}$ starts to deteriorate around the seventh mode in this example. In the SI we display the companion figures for other values of $\alpha$ showing that the breakpoint is changing as well when the tail index is different. This is a consequence of what we show in \ref{fig:bulkedge}, i.e. the threshold to use to discriminate between the structural modes and those coming from randomness is a function of the tail index $\alpha$.}

    \label{fig:eigenvectors1}
\end{figure}

We now turn to the eigenvectors for $k>1$. %Seldom can a description of the eigenvectors be given analytically in terms of the actual entries, but we are able to do so making use of the very same ansatz used for the eigenvalues. 
Plugging $\omega_k$ into Eq.~\eqref{eq:eigenf>}, we get
\begin{equation}
\nu_k(x) = \frac{x^{\alpha/2}}{\alpha}\Re\left[\left(\frac{\alpha}{2}-\i\omega_k\right)\Gamma\left(1-\frac{\alpha}{2}-\i\omega_k\right)x^{\i\omega_k}\right].
\label{eq:eigenf3}
\end{equation}
%\begin{equation}
%\nu_k(x) = -\frac{x^{\alpha/2}}{\alpha}\left(\frac{\alpha^2}{4}+{\omega_k^2}\right)\Re\left[\Gamma\left(-\frac{\alpha}{2}-\i\omega_k\right)x^{\i\omega_k}\right]
%\label{eq:eigenf3}
%\end{equation}
Using Eq.~\eqref{eq_mapping}, we re-express the eigenvector entries as \begin{equation}
    v_k^{(j)}= \sqrt{\frac{n}{j \alpha^2}}
    ~\Re\left[\left(\frac{\alpha}{2}-\i\omega_k\right)\Gamma\left(1-\frac{\alpha}{2}-\i\omega_k\right)\left(\frac{n}{j}\right)^{\frac{\i \omega_k}{\alpha}}\right].\label{eq:eigenvec}
\end{equation}
Since $x^{\i\omega_k}=e^{\i\omega_k\ln x}$ and $j^{\i\omega_k/\alpha}=e^{\i(\omega_k/\alpha)\ln j}$, the structure of the above eigenfunctions and eigenvectors is oscillating, log-periodic in their arguments ($x$ and $j$ respectively), and modulated by a power-law profile ($x^{\alpha/2}$ and $j^{-1/2}$ respectively).
In Fig.~\ref{fig:eigenvectors1} and in SI we show that our analytically derived log-periodic formula tracks the actual eigenvectors of both $\mathbf{P}$ and $\mathbf{A}$ remarkably well for small $k$. 
These eigenvectors do not display traditional delta-like localization, as they are broadly supported on a wide range of weight percentiles with largest amplitude around the hub node ($j=1$), followed by power-law decay $j^{-1/2}$ (see also SI). 
As $k$ grows large enough (below and in SI we provide evidence that this happens for $k$ growing with $n$ but of order at most $\ln n$), we find a departure of the eigenvectors of $\mathbf{A}$ from those of $\mathbf{P}$, due to the entrance of the associated eigenvalues into the random bulk where eigenvectors are expected to be more uniformly spread and less structured. 
%A more detailed characterization of these noisy eigenvectors is beyond the scope of our current analytic approach.

For the parameters of Fig.~\ref{fig:eigenvectors1}, $n=10^4$, and hence $\ln n=\ln(10^4)\simeq 9.2$. The observed deterioration around $k\simeq 7$ is therefore compatible with a cutoff of the form
\[
k_{\mathrm{out}}(n,\alpha)\simeq c(\alpha)\ln n,
\]
with a prefactor $c(\alpha)<1$. More explicitly, if the edge of the bulk generated by $\mathbf H$ is $C_H(\alpha)\sqrt n$, then the last visible outlier is estimated by
\[
k_{\mathrm{out}}(n,\alpha)
=
\max\left\{
k:\alpha\left|\Gamma\left(-\frac{\alpha}{2}+i\omega_k\right)\right|
>
C_H(\alpha)
\right\}.
\]
Since $\omega_k$ is approximately proportional to $k/\ln n$ in the range considered in Eq.~\eqref{sollazzo}, this comparison gives a logarithmic cutoff with an $\alpha$-dependent prefactor. The companion figures in the SI for other values of $\alpha$ show that this breakpoint changes with the tail index, as expected from the dependence of the bulk edge on $\alpha$.

\subsection*{Bulk of the spectrum}
\begin{figure}[t]
    \centering
    \includegraphics[width=\linewidth]{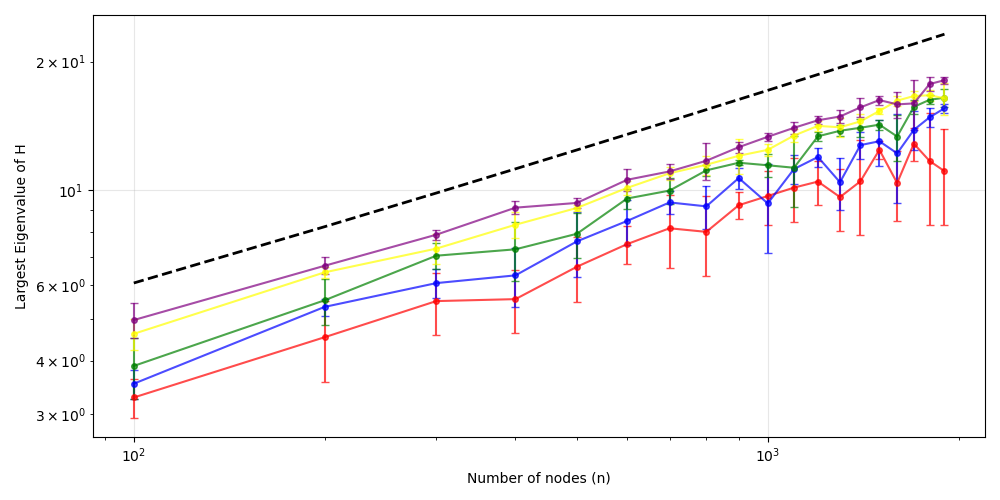}
    \caption{\textbf{Spectral norm of the noise matrix $\mathbf{H}$.}
The largest eigenvalue of $\mathbf{H}=\mathbf{A}-\mathbf{P}$ (i.e.\ the edge of the bulk), averaged over $10$ independent realisations, is plotted as a function of the number of nodes $n$ for five values of the tail index: $\alpha=0.3$ (red circles), $\alpha=0.4$ (blue circles), $\alpha=0.5$ (green circles), $\alpha=0.6$ (yellow circles), and $\alpha=0.7$ (purple circles). The black dashed line shows the analytic upper bound $\sqrt{n}/2$ derived in Eq.~\eqref{eq:upperbound}. All curves grow as $\sqrt{n}$, and the upper bound provides a good envelope across all tested values of $\alpha$, although the $\alpha$-dependent prefactor requires further refinement.}
    \label{fig:bulkedge}
\end{figure}

As we have seen in Fig. \ref{fig:fullspectrum} for the eigenvalues and from Fig. \ref{fig:eigenvectors1} for the eigenvectors (plus more illustrations in SI), there is a value of $k$ beyond which the spectral properties of $\mathbf{P}$ are no longer able to describe those of $\mathbf{A}$.  We refer to this as the `entrance' into the random bulk, i.e. the  noisy part of the spectrum generated primarily by the fluctuation of $\mathbf{H}= \mathbf{A}-\mathbf{P}$.
The analysis of this part of the spectrum is highly non-trivial, and is beyond the scope of this paper. Nevertheless, in SI we show how an analytic characterization of this region is possible using RMT methods. 
In particular, we show that the classical form of the Dyson self-consistency equation for the resolvent is still in place and the usual cavity method yields the same structure of a non-standard method based on Poisson point processes and representation of stable laws, which is well suited for our special regime of weights. 
What we are mainly interested in here is the estimation of where the signal `meets' the bulk. 
In particular, we look for an upper bound for the spectral norm of $\mathbf{H}$, whose elements $H_{ij} = A_{ij}-P_{ij}$ are independent but not identically distributed random variables in the range $[-1,1]$, the variance of $H_{ij}$ being $p_{ij}(1-p_{ij})$. 
The study of such so-called random matrices with a variance profile, or structured random matrices, is of great interest in the mathematical ecology community, and has fostered a lively discussion in the mathematical literature \cite{van2017spectral, latala2005some, schuett2013expectation, cheliotis2025limit, benaych2020spectral, benaych2019largest, bandeira2016sharp, van2017structured}. In particular, we are interested in the bound derived in \cite{bandeira2016sharp}, that can be extended to matrices with entries having a different subgaussian distributions, like the centered Bernoulli of our case. In SI we show that for our case the bound for the spectral norm of $\mathbf{H}$ can be attained as
\begin{equation}
    ||\mathbf{H}|| \lesssim \frac{\sqrt{n}}{2} + \frac{1}{4}\sqrt{\ln{n}} \simeq \frac{\sqrt{n}}{2}\label{eq:upperbound}.
\end{equation}
In Fig.~\ref{fig:bulkedge} we show that the square root growth (that is commonplace in Wigner matrices) seems to be the right one, and that the coefficient coming from the maximization of the row/column sum of the variance profile matrix is already a quite good envelope for all values of $\alpha$. Nonetheless this crude bound needs refinement, at least to include the dependence on $\alpha$. This calls for the full characterization of the bulk spectrum, which we leave to future work.
What is important for our purposes here is noticing that the edge of the bulk lives on the same scale $\sqrt{n}$ as the leading eigenvalues given by Eq.~\eqref{eq:eigenvk}. Since the prefactor of $\sqrt{n}$ for those outliers decays exponentially with $\omega_k$ (and, via Eq.~\eqref{sollazzo}, with $k/\ln n$ -- see SI), we see that the bulk will necessarily `eat out' the leading eigenvalues for some value of $k$ of order at most $\ln n$.

\section*{Discussion} 
In this paper, we provided a characterization of the spectral properties of the Multi-Scale Model (MSM) in the infinite-mean regime ($0<\alpha<1$), which has been recently identified as an annealed random graph model that is exactly invariant under node aggregation~\cite{garuccio2023multiscale,avena2022inhomogeneous}. This invariance is a key concept in the field of network renormalization~\cite{gabrielli2025network}, but its connection with discrete scale invariance and the associated idea of complex scaling exponents~\cite{logperiodic1,logperiodic2} had not been investigated so far. 
We showed that the leading eigenvalues of the expected adjacency matrix $\mathbf{P}$ scale as $\sqrt{n}$ and accurately proxy the outlier eigenvalues (whose number is growing and of order at most $\ln n$) of the realized adjacency matrix $\mathbf{A}$. 
Remarkably, these outliers have an approximately constant spacing on a logarithmic scale of their index $k$ and can be identified as the intersection between the real axis and a \emph{spira mirabilis} in the complex plane, a self-similar object naturally emerging in presence of complex scaling exponents and stable laws with complex stability index~\cite{alexeev2022stable,alexeev2023stable}. 
We also derived a closed-form expression for the outlier eigenvectors, showing they are delocalized but biased toward high-weight nodes. The leading eigenvectors beyond the principal (Perron–Frobenius) one exhibit log-periodic oscillations,
a feature naturally produced by discrete scale invariance~\cite{logperiodic1,logperiodic2,logperiodic3} which is here present in the property that the same random graph model describes exactly (up to a parameter rescaling) any coarse-grained version of the original graph, where sets of nodes of the same cardinality are aggregated into supernodes~\cite{garuccio2023multiscale,avena2022inhomogeneous}.
Our results qualitatively confirm the findings of~\cite{nadakuditi2013spectra, benaych2011eigenvalues}, while partially extending them to the regime of infinite-rank perturbations of random matrices as in more recent attempts~\cite{afanasiev2025asymptotic, huang2018mesoscopic}.

The leading eigenvalues therefore have a direct structural interpretation: they identify the deterministic multi-scale component that survives edge randomness. Estimating them accurately is useful for detecting discrete scale invariance, for distinguishing structural outliers from the random bulk, and potentially for inferring the tail index or renormalization parameters from spectral data. The same strategy may also be useful for other heavy-tailed inhomogeneous random graphs, where the expected adjacency matrix is not finite-rank but has a hierarchy of slowly decaying structural modes.
In summary, the MSM illustrates how infinite-mean heterogeneity generates spectral phenomena far from classical paradigms. 
%Our analysis suggests parallels with other models such as Chung–Lu, where broad degree distributions induce hub-driven outliers. 
Open directions include testing the universality of the $\sqrt{n}$ scaling and log-periodicity across other networks and models. 
More broadly, our results provide a framework for understanding spectra in multi-scale and heavy-tailed systems, and lay the ground for future work in these fields.

\section*{Methods}

\subsection*{Graph model and parameter choices}

Numerical simulations are performed on the annealed Multi-Scale Model (MSM) defined in the Results section. Each graph on $n$ nodes is generated as follows. Node fitnesses $\{x_i\}_{i=1}^n$ are drawn independently from the Pareto($\alpha$) distribution $\rho_\alpha(x)=\alpha x^{-1-\alpha}$ for $x\ge1$, with $0<\alpha<1$. The connection probability between nodes $i$ and $j$ is $p_{ij}=1-\exp(-\varepsilon_n x_i x_j)$ with $\varepsilon_n=n^{-1/\alpha}$. Each edge $(i,j)$ for $i<j$ is then drawn independently as a Bernoulli($p_{ij}$) variable, and the resulting symmetric adjacency matrix $\mathbf{A}$ is stored in sparse format. In addition to random fitness assignments, we also consider deterministic fitness assignments $x_i=(n/i)^{1/\alpha}$ (nodes ranked in decreasing order) to suppress sample-to-sample fluctuations at finite $n$.

\subsection*{Eigenvalue and eigenvector computations}

The expected adjacency matrix $\mathbf{P}$ with entries $P_{ij}=p_{ij}$ is stored as a dense matrix for moderate $n$ ($n\le 10^4$). Leading eigenvalues and eigenvectors of both $\mathbf{A}$ and $\mathbf{P}$ are computed using standard symmetric eigensolver routines (scipy.linalg.eigh). The noise matrix $\mathbf{H}=\mathbf{A}-\mathbf{P}$ is assembled and its spectral norm (largest eigenvalue) is computed over multiple independent realisations to estimate the bulk edge. Results shown in the figures use $n=10^4$ and are averaged over 5--10 independent realisations as indicated in the figure captions.

\subsection*{Analytic predictions}

The analytic eigenvalue predictions (Eq.~\eqref{eq:eigenvk}) and the logarithmic spiral (Eq.~\eqref{eq:spiral}) are evaluated numerically by solving the implicit equation (Eq.~\eqref{eq:admissible}) for $\omega_k$ using the secant root-finding algorithm applied to the argument of the Gamma function. The approximate linear formula (Eq.~\eqref{sollazzo}) is used as an initial guess. Eigenvector predictions (Eq.~\eqref{eq:eigenvec}) are evaluated at the node ranks $j=1,\dots,n$ using the same $\omega_k$ values. All special-function evaluations (Gamma function, its argument) are performed using standard scientific computing libraries in python.

\backmatter

\bmhead{Supplementary information}

We provide a Supplementary Information (SI) file in which we show in detail all the derivations that for the sake of clarity and brevity could not fit the main body of the paper.

\bmhead{Data Availability}

The datasets generated and/or analysed during the current study are not publicly available due to the fact that the numerical simulations can be easily reproduced using the model definition provided in the Methods section, but are available from the corresponding author on reasonable request.

\bmhead{Code Availability}

The code used to generate the numerical simulations and figures in this paper is available from the corresponding author upon reasonable request.

\bmhead{Acknowledgements}
This publication is part of the project ``Redefining renormalization for complex networks'' with file number OCENW.M.24.039 of the research programme Open Competition Domain Science Package 24-1, which is (partly) financed by the Dutch Research Council (NWO) under the grant https://doi.org/10.61686/PBSEC42210.
This work is also supported by the European Union - NextGenerationEU and funded by the Italian Ministry of University and Research (MUR) - National Recovery and Resilience Plan (Piano Nazionale di Ripresa e Resilienza, PNRR), projects ``Strengthening the Italian RI for Social Mining and Big Data Analytics'' (SoBigData.it, grant IR0000013, DN. 3264, 28/12/2021), ``Reconstruction, Resilience and Recovery of Socio-Economic Networks'' (RECON-NET, EP\_FAIR\_005 - PE0000013 ``FAIR'' - PNRR M4C2 Investment 1.3), ``A Multiscale integrated approach to the study of the nervous system in health and disease'' (MNESYS, PE0000006, DN. 1553, 11/10/2022).

\bmhead{Author Contributions}

A.C., R.S.H and D.G\ conceived and performed the analytical derivations and prepared the manuscript. A.C. conducted the numerical simulations, R.S.H.\ contributed to the mathematical framework, D.G.\ supervised the project and critically revised the manuscript. All authors read and approved the final manuscript.

\bmhead{Competing Interests}

The authors declare no competing financial or non-financial interests.


\bibliography{sn-bibliography}% common bib file

%% ============================================================
%% Figure Legends
%% ============================================================

%%%%%%%%%%%%%%%%%%%%%%%%%%%%%%%%%%%%%%%%%%%%%%%%%%%%%%%%%%%%%%%
%
%    SUPPLEMENTARY MATERIAL
%
%%%%%%%%%%%%%%%%%%%%%%%%%%%%%%%%%%%%%%%%%%%%%%%%%%%%%%%%%%%%%%%
%%%%%%%%%%%%%%%%%
\clearpage
\newpage
%\appendix
\setcounter{equation}{0}

\setcounter{figure}{0}
\setcounter{table}{0}
\setcounter{page}{1}
\setcounter{section}{0}
%\makeatletter

\renewcommand{\theequation}{S\arabic{equation}}
\renewcommand{\thetable}{S\arabic{table}}
\renewcommand{\thefigure}{S\arabic{figure}}
\renewcommand{\thesection}{S.\Roman{section}} 
\renewcommand{\thesubsection}{\thesection.\roman{subsection}}

%\vspace{\columnsep}
{\center
\textbf{SUPPLEMENTARY INFORMATION}\\
$\quad$\\
accompanying the paper\\
\emph{``Spectra of random graphs with discrete scale invariance''}\\
by A. Catanzaro, R.S. Hazra and D. Garlaschelli\\
$\quad$\\
$\quad$\\
}
%\twocolumngrid

\section*{Asymptotics of the Laplace Transform}
\label{app:LT}

We begin by explicitly computing the action of our integral operator on the function $\nu_1(x)$:
\begin{eqnarray}
\mathbf{\Pi}_{\alpha,n}[\nu_1(x)]&=&1-\varphi_{\beta}(\varepsilon_n x)\nonumber\\
&=&1-\int_1^{+\infty}dy~\rho_\beta(y)~ e^{-\varepsilon_n x y}\nonumber\\
&=&1-\beta \int_1^{+\infty}dy~y^{-1-\beta}~ e^{-\varepsilon_n xy}\nonumber\\
&=&1-\beta \int_{\varepsilon_n x}^{+\infty}dt~{\left(\frac{t}{\varepsilon_n x}\right)}^{-1-\beta}~ \frac{e^{-t}}{\varepsilon_n x}\nonumber\\
&=&1-\frac{\beta}{(\varepsilon_n x)^{-\beta}} \int_{\varepsilon_n x}^{+\infty}dt~{t}^{-1-\beta}~ {e^{-t}}\nonumber\\
&=&1-\frac{\beta}{(\varepsilon_n x)^{-\beta}}~ \Gamma(-\beta,\varepsilon_n x),
\label{eq:prod1}
\end{eqnarray}
where we have made the change of variables $\varepsilon_n x y\to t$ and introduced the incomplete Gamma function $\Gamma(z,s)$, which has the following definition and asymptotic relation with the (complete) Gamma function $\Gamma(z)\equiv\Gamma(z,0)$ for $\Re[z]<0$: 
\begin{equation}
\Gamma(z,s)\equiv\int_s^{+\infty}dt~t^{z-1}e^{-t}~\sim~\Gamma(z)+\frac{s^{z}}{z},\qquad s\to0^+.
\label{Seq:scaling}
\end{equation}
Since in our case $\Re[z]=-\alpha/2<0$, the use of the above relation is justified.
Plugging Eq.~\eqref{Seq:scaling} into Eq.~\eqref{eq:prod1}, we obtain for $\varepsilon_n x\to 0^+$:
\begin{eqnarray}
   \mathbf{\Pi}_{\alpha,n}[\nu_1(x)] &\sim&
   1-\frac{\beta}{(\varepsilon_n x)^{-\beta}}\left[\Gamma(-\beta)-\frac{(\varepsilon_n x)^{-\beta}}{\beta}\right]\nonumber\\
   &=&-\beta~(\varepsilon_n x)^{\beta}~\Gamma(-\beta)\nonumber\\
   &=&(\varepsilon_n x)^{\beta}~\Gamma(1-\beta)
   \label{eq:prod2}
\end{eqnarray}
which is the relation used in the main text in the case $k=1$.

\section*{Eigenvalues for $k>1$}
We will use the fact that the Gamma function of the complex conjugate of its argument is the complex conjugate of the Gamma function. 
The discussion in the previous section extends to complex-valued parameters $\beta^\pm_k=\gamma_k\pm\i\omega_k$, leading to the following generalization of Eq.~\eqref{eq:prod2}: \begin{equation}
\mathbf{\Pi}_{\alpha,n}[\nu^\pm_k(x)]\sim (\varepsilon_nx)^{\beta_k^\pm} ~\Gamma(1-\beta_k^\pm).
%=(\varepsilon_nx)^{\gamma}(\varepsilon_nx)^{\pm i\omega} ~\Gamma(1-\gamma\mp i\omega_k^\pm)
      \label{eq:test}
\end{equation} 
Given the candidate eigenfunction $\nu_k(x)=\Re\left[\Gamma(1-\beta_k^-)~\nu^+_k(x)\right]$, Eq.~\eqref{eq:test} leads to the following asymptotics:
\begin{equation}
\mathbf{\Pi}_{\alpha,n}[\nu_k(x)]\sim \frac{1}{2}~\Gamma(1-\beta_k^+)~\Gamma(1-\beta_k^-)~\left[(\varepsilon_nx)^{\beta_k^+} +(\varepsilon_nx)^{\beta_k^-}\right]=\left|\Gamma(1-\beta_k^+)\right|^2~\Re\left[(\varepsilon_nx)^{\beta_k^+} \right].%\\
%&=&|\Gamma(1-\gamma_k\pm i\omega_k)|^2(\varepsilon_nx)^{\gamma_k}\left[(\varepsilon_nx)^{+ i\omega}+(\varepsilon_nx)^{- i\omega}\right],
      \label{eq:lhs}
\end{equation} 
To enforce the eigenvector equation $\mathbf{\Pi}_{\alpha,n}[\nu_k(x)]=\frac{\lambda_k}{n}~\nu_k(x)$, the above expression has to equal 

\begin{align}
    \frac{\lambda_k}{n}~\nu_k(x) &=
   \frac{\lambda_k}{2n\alpha}\left[\Gamma(1-\beta_k^-)~
    {\beta_k^+}~x^{\alpha-\beta_k^+}   
    +\Gamma(1-\beta_k^+){\beta_k^-}~x^{\alpha-\beta_k^-}    
    \right] 
    \\
    &= \frac{\lambda_k}{n\alpha}\Re\left[\Gamma(1-\beta_k^-)~
    {\beta_k^+}~x^{\alpha-\beta_k^+}   
    \right].
    \label{eq:rhs}
\end{align}

\begin{figure}[]
    \centering 
    \includegraphics[width=0.6\linewidth]{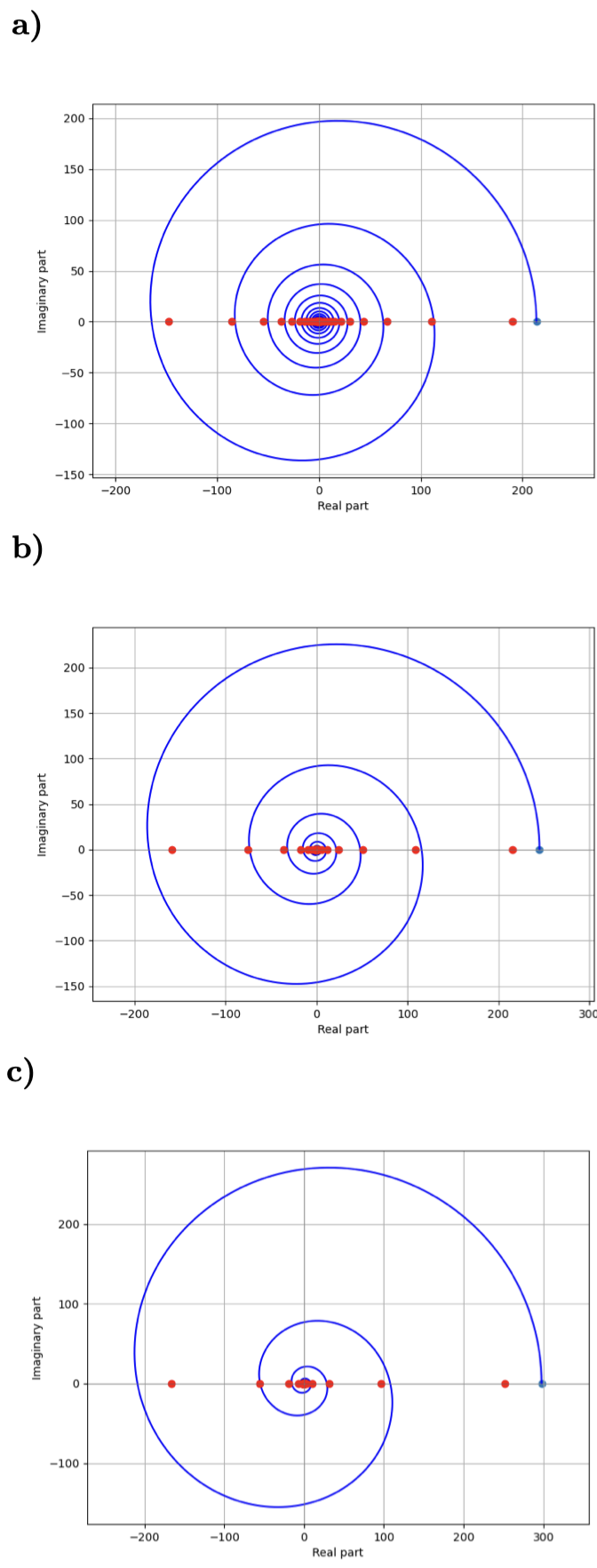}
    \caption{\textbf{Eigenvalues of the expected adjacency matrix $\mathbf{P}$ lie at the intersections of the real axis with a logarithmic \emph{spira mirabilis} in the complex plane.} \textbf{a)} $\alpha = 0.2$. \textbf{b)} $\alpha = 0.5$. \textbf{c)} $\alpha = 0.8$.
    The actual real-valued outlier eigenvalues of $\mathbf{P}$ (red dots) are shown together with the logarithmic spiral $\sigma^+(\omega)$ (blue curve) defined in Eq.~\eqref{eq:spiral}, plotted parametrically in the complex plane as a function of $\omega\ge0$. Here $n=10^4$ and $\alpha=0.2, 0.5, 0.8$ from top to bottom.}
    \label{fig:spirals}
\end{figure}

This requires either $\beta_k^+\equiv \alpha-\beta_k^+$ and $\beta_k^-\equiv \alpha-\beta_k^-$, which would however retrieve the case $k=1$ that we have already considered for the largest eigenvalue ($\beta_k^\pm=\gamma_k=\alpha/2$ and $\omega_k=0$, since $\alpha\in\mathbb{R}$), or $\beta_k^+\equiv \alpha-\beta_k^-$ and $\beta_k^-\equiv \alpha-\beta_k^+$, which still leads to $\gamma_k=\alpha/2$ (for all $k$) but puts no restriction on $\omega_k$, so that $\beta_k^\pm=\alpha/2\pm \i\omega_k$.
The defining restriction on $\omega_k$ is obtained by matching the other terms in Eqs.~\eqref{eq:lhs} and~\eqref{eq:rhs}, yielding the two simultaneous conditions \begin{equation}
    \lambda_k = -{\alpha}~ {\Gamma\left(-\frac{\alpha}{2}\mp\i\omega_k\right)}~n^{1/2~\pm~\i\omega_k/\alpha}= -{\alpha}~ {\Gamma\left(-\frac{\alpha}{2}\mp\i\omega_k\right)}~n^{1/2}~ e^{\pm\i(\omega_k/\alpha)\ln n}.\label{eq:lambda12}
\end{equation} 
Since the right-hand sides of the above expression are complex conjugates of each other, the two simultaneous  conditions ensure that $\lambda_k$ is real, as required for the eigenvalues of a symmetric matrix. 
This means that $\omega_k$ has to be such that $\arg{[\lambda_k]}=0,\pi$ modulo $2\pi$, i.e.
\begin{equation}
\arg{\left[{\Gamma\left(-\frac{\alpha}{2}+\i\omega_k\right)}\right]}=
    \frac{\omega_k}{\alpha}\ln n-k\pi, \quad k\in\mathbb{N},\qquad -\pi\le\arg{\left[{\Gamma\left(-\frac{\alpha}{2}+\i\omega_k\right)}\right]}<+\pi,
\label{Seq:admissible}
\end{equation}
which coincides with Eq.~\eqref{eq:admissible} in the main text.
Note that the (necessarily real-valued) eigenvalues $\lambda_k$ ($k\ge1$) can be found graphically as the intersections between the real axis and the two  spirals 
\begin{equation}
\sigma^\pm(\omega)\equiv-{\alpha}~ \Gamma\left(-\frac{\alpha}{2}\mp \i\omega\right)n^{1/2}e^{\pm \i(\omega/\alpha)\ln n}    
\label{Seq:spirals}
\end{equation}
defined parametrically in the complex plane by extending the $\omega_k$ in Eq.\eqref{eq:lambda12} to any real value $\omega\ge 0$. 
See Fig.~\ref{fig:spirals} for an illustration of this result, for three representative values of the index $\alpha =0.2, 0.5, 0.8$.
When the solutions for $\omega_k$ that realize Eq.~\eqref{Seq:admissible} (or equivalently the intersections) are inserted into Eq.\eqref{eq:lambda12}, they produce all the eigenvalues as
\begin{equation}  \label{Seq:eigenvk}  
     \lambda_k = (-1)^k~\alpha \left| {\Gamma\left(-\frac{\alpha}{2}+\i\omega_k\right)}\right|~ n^{1/2},
\end{equation}
proving the expression in the main text, which nicely extends the expression for $\lambda_1$ to values $k> 1$.

\subsection*{Graphical inspection}
\begin{figure}[b]
    \centering
    \includegraphics[width=0.8\linewidth]{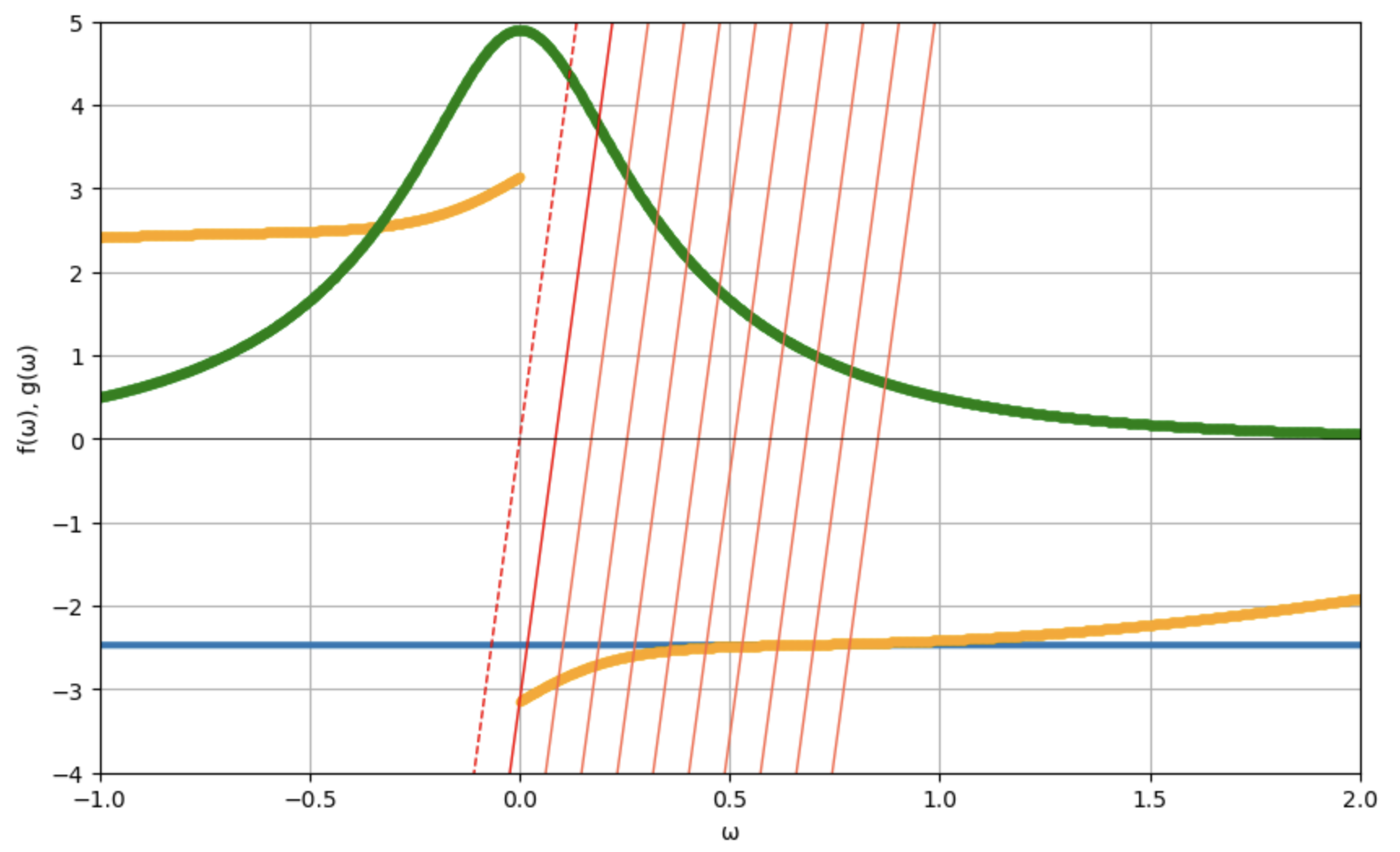}
    \caption{\textbf{Graphical Solution of Eq.~\eqref{Seq:admissible}}. In orange, the function $f(\omega)$ defined in  Eq.~\eqref{Seq:f_omega}, while the straight lines $g_k(\omega)$ defined by Eq.~\eqref{Seq:g_omega} are in red, with $k$ from $0$ to $10$. The dashed red line corresponds to $g_0(\omega)$, which has no intersection $\omega_0$ with $f(\omega)$. The solid red line at its immediate right is $g_1(\omega)$, which intersects $f(\omega)$ at the  exact value $\omega_1 =0$ where  $f(\omega_1)=g_1(\omega_1)=-\pi$. Further on the right, the next straight lines $g_k(\omega)$ (with $k>1$) are illustrated in light red.
    The solution to Eq.~\eqref{Seq:admissible} for these values of $k$ can be obtained graphically where the orange and light red lines meet. In this regime, $f(\omega)$ can be approximated as the  horizontal line $f(\omega)\approx\phi_\alpha\equiv f(\omega_\alpha)$ (blue line), where the value $\omega_\alpha$ is the stationary point such that $f'(\omega_\alpha)=0$. 
    The horizontal line is no longer a good approximation for $f(\omega_k)$ for larger $k$, but that is already the regime for which the value of $\left|{\Gamma\left(-\frac{\alpha}{2}+\i\omega_k\right)}\right|$ (plotted in green) is exponentially suppressed, as one can verify from Eqs.~\eqref{Seq:expodacay}.}
    \label{fig:mecojoni}
\end{figure}

Unfortunately, Eq.~\eqref{Seq:admissible} cannot be solved explicitly for $\omega_k$, except for $k=0$ (for which there is no solution for $\omega_0$) and $k=1$, for which $\omega_1=0$ is a solution, leading indeed to the value $\lambda_1$ when inserted into Eq.\eqref{Seq:eigenvk}: indeed,  ${\Gamma\left(-\frac{\alpha}{2}\right)}<0$, so that $\arg{\Gamma\left(-\frac{\alpha}{2}\right)}=-\pi$ solves Eq.~\eqref{Seq:admissible} for $k=1$ and $\omega_1=0$.
The solutions for $k>1$ can be inspected graphically, either as the aforementioned intersections between the real axis and the complex spiral, or as the intersections between the two curves defined as the l.h.s. and r.h.s. of Eq.~\eqref{Seq:admissible}.
Specifically, the l.h.s. defines the function
\begin{equation} \label{Seq:f_omega}  f(\omega)\equiv\arg{{\Gamma\left(-\frac{\alpha}{2}+\i\omega\right)}},
\end{equation}
while the r.h.s. defines the family of straight lines \begin{equation}\label{Seq:g_omega}  
g_k(\omega)\equiv{\omega}\frac{\ln n}{\alpha}-k\pi,\quad k\ge0.
\end{equation}
The solutions $\omega_k$ to Eq.~\eqref{Seq:admissible} are the values at the intersections
\begin{equation}
    f(\omega_k)=g_k(\omega_k),\quad k\ge0.
    \label{Seq:intersection}
\end{equation}
In Fig.\ref{fig:mecojoni} we plot $f(\omega)$ and $g_k(\omega)$ for a few small values of $k$.
The plot confirms that there is no intersection point $\omega_0$ for $k=0$, while the intersection for $k=1$  is found at $\omega_1=0$, for which $f(\omega_1)=g_1(\omega_1)=-\pi$.
What is more useful from the graphical inspection is that the next few intersections for $1<k\le k^*$ occur at a range of values $\omega_k$ close to some value $\omega_\alpha$ at which $f(\omega)$ has a stationary point and varies very slowly around the corresponding value $\phi_\alpha\equiv f(\omega_\alpha)$, which is approximately a local `plateau' for the function (the value $k^*$ is defined such that for $k>k^*$ the value $f(\omega_k)$ is away from the plateau).
We can therefore approximate the solutions $\omega_k$ for the first few values $k=2,3,\dots,k^*$ as the points at which $g_k(\omega)$ intersects the constant line $f(\omega)\approx \phi_\alpha$, 
and replace Eq.~\eqref{Seq:intersection} with
\begin{equation}
    \phi_\alpha\approx g_k(\omega_k),\quad 1<k\le k^*.
    \label{Seq:approxintersection}
\end{equation}
which, using Eq.~\eqref{Seq:g_omega}, is solved by the values
\begin{equation}
    \omega_k \approx \alpha\frac{k\pi +\phi_\alpha}{\ln n},\quad \phi_\alpha\equiv f(\omega_\alpha),\quad 1<k\le k^*.
    \label{Seq:sollazzo}
\end{equation}
To complete our characterization of these approximate solutions, we need to estimate both $k^*$ and $\omega_\alpha$.

\subsection*{Estimation of $\omega_\alpha$}
\begin{figure}[h]
    \centering
    \includegraphics[width=0.5\linewidth]{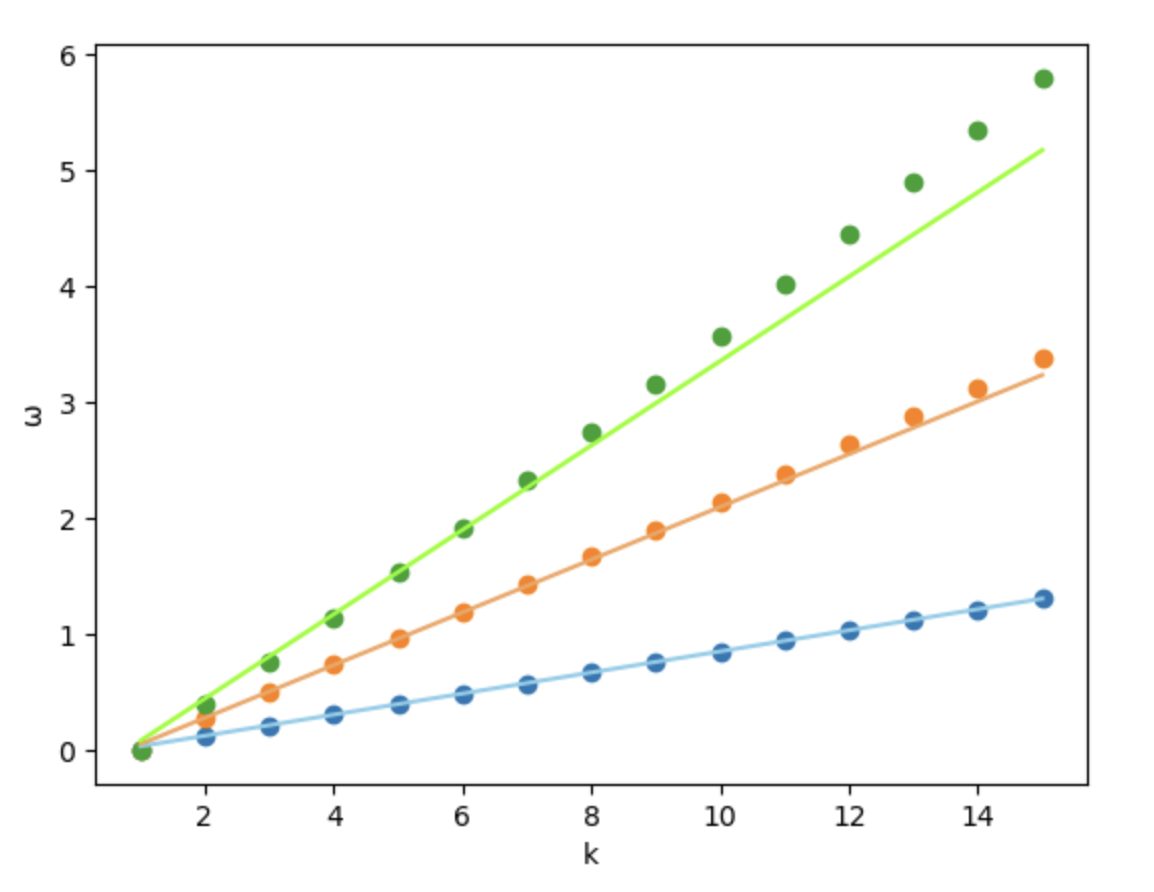}
    \caption{\textbf{Values of $\omega_k$  as a function of $k$.} The isolated points are the actual, numerical solutions of Eq.~\eqref{Seq:admissible}, for $\alpha = 0.2, 0.5, 0.8$ (in blue, orange and green respectively). The solid lines are the approximate linear solutions defined by  Eq.~\eqref{Seq:sollazzo} with $k>1$ (plus the exact solution $\omega_1=0$) for the same  values of $\alpha =0.2, 0.5, 0.8$ (and same colors).}
    \label{fig:omega_vs_k}
\end{figure}

We start from the determination of $\omega_\alpha$.
Since the latter is defined as a stationary point of $f(\omega)$, we look for its value by calculating the derivative $f'(\omega)$ and requiring $f'(\omega_\alpha)=0$, or more precisely:
\begin{equation} \label{Seq:derivative}  \frac{\partial}{\partial\omega}\arg{{\Gamma\left(-\frac{\alpha}{2}+\i\omega\right)}}\Big|_{\omega=\omega_\alpha}=0.
\end{equation}
We now rearrange the identity $z=|z|e^{\i\arg{z}}=\sqrt{z~\bar{z}}~e^{\i\arg{z}}$ (valid for every complex number $z$, where $\bar{z}$ is the complex conjugate of $z$) to express 
\begin{equation}
    \arg{z}=-\i\ln\frac{z}{|z|}=\frac{\i}{2}\ln\frac{\bar{z}}{z}.
    \label{Seq:logarg}
\end{equation}
If $z=\arg{{\Gamma\left(-\frac{\alpha}{2}+\i\omega\right)}}$, we get
\begin{equation}
    \arg{{\Gamma\left(-\frac{\alpha}{2}+\i\omega\right)}}=\frac{\i}{2}\ln\frac{\Gamma\left(-\frac{\alpha}{2}-\i\omega\right)}{\Gamma\left(-\frac{\alpha}{2}+\i\omega\right)}=\frac{\i}{2}\ln{\Gamma\left(-\frac{\alpha}{2}-\i\omega\right)}-\frac{\i}{2}\ln{\Gamma\left(-\frac{\alpha}{2}+\i\omega\right)}.
    \label{Seq:argGamma}
\end{equation}
In order to be able to differentiate the above expression with respect to $\omega$, we need a closed-form expression for $\ln\Gamma(z)$ (also called the `log-Gamma function'), which is available as follows:
\begin{equation}
    \ln\Gamma(z)=-\gamma z-\ln z+\sum_{m=1}^{+\infty}\left[\frac{z}{m}-\ln\left(1+\frac{z}{m}\right)\right],
    \label{Seq:log-gamma}
\end{equation}
where
\begin{equation}
\gamma\equiv\lim_{m\to+\infty}\left[-\ln m+\sum_{l=1}^{m}\frac{1}{l}\right]= 0.5772\dots    
\end{equation}
is the Euler-Mascheroni constant. 
Using Eq.~\eqref{Seq:log-gamma} with $z=-\alpha/2\pm\i\omega$ into Eq.~\eqref{Seq:argGamma}, we arrive at
\begin{equation}
    \arg{{\Gamma\left(-\frac{\alpha}{2}+\i\omega\right)}}=
    -\gamma \omega-\arg{{\left(-\frac{\alpha}{2}+\i\omega\right)}}+\sum_{m=1}^{+\infty}\left[\frac{\omega}{m}-\arg{{\left(m-\frac{\alpha}{2}+\i\omega\right)}}\right].
    \label{Seq:ciotola}
\end{equation}
Now, in order to differentiate with respect to $\omega$ as prescribed by Eq.~\eqref{Seq:derivative}, we preliminarily calculate 
\begin{equation} \label{Seq:derivarg}  \frac{\partial}{\partial\omega}\arg{{\left(m-\frac{\alpha}{2}+\i\omega\right)}}=\frac{\i}{2}\frac{\partial}{\partial\omega}\ln\left(\frac{m-\frac{\alpha}{2}-\i\omega}{m-\frac{\alpha}{2}+\i\omega}\right)=\frac{m-\frac{\alpha}{2}}{\left(m-\frac{\alpha}{2}\right)^2+\omega^2},
\end{equation}
where we have used Eq.~\eqref{Seq:logarg} again.
Using the above expression into Eq.~\eqref{Seq:ciotola}, we obtain
\begin{equation}
    \frac{\partial}{\partial\omega}\arg{{\Gamma\left(-\frac{\alpha}{2}+\i\omega\right)}}=
    -\gamma +\frac{\frac{\alpha}{2}}{\left(\frac{\alpha}{2}\right)^2+\omega^2}
    +\sum_{m=1}^{+\infty}\left[\frac{1}{m}-\frac{m-\frac{\alpha}{2}}{\left(m-\frac{\alpha}{2}\right)^2+\omega^2}\right].
    \label{Seq:ciotolona}
\end{equation}
Now, imposing Eq.~\eqref{Seq:derivative} translates into the following: 
\begin{equation}
-\gamma +\frac{\frac{\alpha}{2}}{\left(\frac{\alpha}{2}\right)^2+\omega_\alpha^2}
    +\sum_{m=1}^{+\infty}\left[\frac{1}{m}-\frac{m-\frac{\alpha}{2}}{\left(m-\frac{\alpha}{2}\right)^2+\omega_\alpha^2}\right]=0,
\end{equation}
which, if the sum is truncated to any finite order for $m$, can be solved algebraically for $\omega_\alpha$.
In general, the number of solutions of the above equation (i.e. stationary points of the function $f$) may depend on $\alpha$. How many (and how closely) such solutions are recovered after truncating the sum will depend on the order of truncation.
For our purposes here (i.e. identifying the first stationary point as shown in Fig.~\ref{fig:mecojoni}), we simply truncate to zeroth order and neglect the sum altogether, to get the condition 
\begin{equation}
-\gamma +\frac{\frac{\alpha}{2}}{\left(\frac{\alpha}{2}\right)^2+\omega_\alpha^2}=0,
\end{equation}
which is solved by the single (positive) value
\begin{equation}
\omega_\alpha\equiv\sqrt{\frac{\alpha}{2}\left(\frac{1}{\gamma}-\frac{\alpha}{2}\right)}.
\end{equation}
This finally allows us to compute
\begin{equation}
    \phi_\alpha\equiv f(\omega_\alpha)=f\left(\sqrt{\frac{\alpha}{2}\left(\frac{1}{\gamma}-\frac{\alpha}{2}\right)}\right),
\end{equation}
which, when inserted into Eq.~\eqref{Seq:sollazzo}, produces our approximate solutions for $\omega_k$ to Eq.~\eqref{Seq:admissible} with $1<k\le k^*$.
A comparison between the actual solutions $\omega_k$ of Eq.~\eqref{Seq:admissible} and the approximate solutions given in \eqref{Seq:sollazzo}
is shown in Fig.~\ref{fig:omega_vs_k}, confirming a very good agreement, especially for smaller values of $\alpha$.
Note that our only following use of these approximate solutions is the approximate proportionality between $\omega_k$ and $k$, which is used below to prove the exponential decay of the prefactor of the leading eigenvalues with $k$.
This linearity only requires that $f(\omega)$ is sufficiently flat around some value $\phi_\alpha$ before $\left| {\Gamma\left(-\frac{\alpha}{2}+\i\omega\right)}\right|$ is exponentially suppressed (see Fig.~\ref{fig:mecojoni}).
A more accurate estimate of $\phi_\alpha$ would improve the estimate of $\phi_\alpha$ used in Eq.~\eqref{Seq:sollazzo} (hence improving the slope of the predicted lines in Fig.~\ref{fig:omega_vs_k}),  but not modify the approximate proportionality between $\omega_k$ and $k$, which is all we need in what follows.

\subsection*{Exponential decay of the prefactor of the leading eigenvalues}
We now turn to the estimation of the order of $k^*$.
Let us consider the Stirling approximation \cite[\href{http://dlmf.nist.gov/5.11.E9}{(5.11.9)}]{NIST:DLMF} for the absolute value of the Gamma function: 
\begin{equation}
    \left| \Gamma(a+ib)\right| \sim \sqrt{2 \pi} |b|^{a-\frac{1}{2}}e^{-\frac{\pi|b|}{2}},
\end{equation}
where $a,b\in\mathbb{R}$. 
Choosing $a=-\alpha/2$ and $b=\omega_k\ge 0$, this implies
\begin{equation}
   \left| {\Gamma\left(-\frac{\alpha}{2}+\i\omega_k\right)}\right| \sim \sqrt{2 \pi}~\omega_k^{-\frac{\alpha+1}{2}}e^{-\frac{\pi\omega_k}{2}},
\end{equation}
showing an exponential decay with increasing $\omega_k$.
This exponential decay is visible in Fig.~\ref{fig:mecojoni}.
If we further use the expression we derived in Eq.~\eqref{Seq:sollazzo} for the approximate solutions $\omega_k$ for $1<k\le k^*$, we see that the above behaviour translates into an exponential decay as a function of $k/\ln n$:
\begin{equation}
   \left| {\Gamma\left(-\frac{\alpha}{2}+\i\omega_k\right)}\right| \sim  \sqrt{2 \pi } \left(\alpha \frac{k \pi + \phi_{\alpha}}{\ln n }\right)^{- \frac{\alpha+1}{2}}  e^{- \frac{\alpha \pi}{2} \frac{k \pi + \phi_{\alpha}}{\ln n}}
,\quad 1<k\le k^*.
\label{eq:lampredotto}
\end{equation}
Looking at Fig.~\ref{fig:mecojoni}, we see that when $k\gtrsim k^*$ (i.e. when the intersections between $f(\omega)$ and $g_k(\omega)$ depart from the approximate plateau of constant height $\phi_\alpha$), the value of $\left| {\Gamma\left(-\frac{\alpha}{2}+\i\omega\right)}\right|$ is already exponentially suppressed.
Therefore Eq.~\eqref{eq:lampredotto} indicates that $k^*$ is of order $\ln n$.
Obviously, the magnitude of the corresponding eigenvalue $\lambda_k$ of $\mathbf{P}$, which is given by Eq.~\eqref{Seq:eigenvk}, is also exponentially (in $k/\ln n$) suppressed:
\begin{eqnarray}
     |\lambda_k| &=& \alpha\left| {\Gamma\left(-\frac{\alpha}{2}+\i\omega_k\right)}\right| \sqrt{n} \\
     &\sim& \alpha~ \sqrt{2 \pi n} \left(\alpha \frac{k \pi + \phi_{\alpha}}{\ln n }\right)^{- \frac{\alpha+1}{2}}  e^{- \frac{\alpha \pi}{2} \frac{k \pi + \phi_{\alpha}}{\ln n}}.\label{Seq:expodacay}
\end{eqnarray}
This means that, for $k=O(\ln n)$, the approximation in Eq.~\eqref{Seq:sollazzo} breaks down, but on the other hand the corresponding values of $\lambda_k$ become exponentially small. 
This is also where, moving from the expected matrix $\mathbf{P}$ to a  realized random matrix $\mathbf{A}$, the randomness will create noisy eigenvalues `eating up' the eigenvalues of $\mathbf{P}$.\\

Indeed, generalized BBP-type results typically assume that the expected matrix $\mathbf{P}$ has a fixed, finite rank in order to guarantee a clean separation between informative (structural) eigenvalues and the random bulk. 
In our model this assumption fails literally, since $\mathbf{P}$ is not finite-rank. Nevertheless, both our numerics and the above arguments indicate that the $k$-th largest (in absolute value) eigenvalue $\lambda_k$ of $\mathbf{P}$ decays rapidly with $k/\ln n$, so that the rank of the matrix is \emph{effectively} of order $\ln n$, and hence $o(n)$. This aligns with recent literature suggesting that a BBP-like transition can persist to order $o(n)$, hence also beyond strict finite rank~\cite{afanasiev2025asymptotic, huang2018mesoscopic}.
Numerics (see main text) indicate that the edge of the random spectrum contributed by $\mathbf{H}$ scales as $C(\alpha)\sqrt{n}$ with $C(\alpha)$ independent of $n$. Consequently, only those signal eigenvalues that are not exponentially suppressed can rise above the bulk, i.e., indices up to $k^*\sim \ln n$. 
In this sense, the signal behaves as if it had \emph{effective rank} $\Theta(\ln n)$: the spectrum of $\mathbf{A}$ consists of order $\ln n$ outliers (delta spikes from $\mathbf{P}$) on top of the bulk determined by $\mathbf{H}$.

\section*{Bulk Spectrum Characterization}\label{app:diago}

In the MSM, because edge probabilities $p_{ij}$ vary across pairs (especially spanning many orders of magnitude due to the weights heterogeneity), the matrix $\mathbf{H}$ does not fall into standard Wigner or Erd\H{o}s-R\'enyi classes. Classical results like the Wigner semicircle law or variants for sparse graphs do not directly apply. Nonetheless, one can attempt to derive a self-consistent equation for the Stieltjes transform. Here we show how the cavity method is able to produce a self-consistent equation for the resolvent of the Spectral Density, and how the derived formula can lead to an equation to estimate the Bulk edge. We also show a different approach based on Poisson point processes, that leads to the very same structure of the self-consistency equation. Finally, we compute a crude bound for the bulk edge and show that the scaling is the one we  wrote in Eq.~\eqref{eq:upperbound}. 

\subsection*{Cavity Approach}

To locate the continuous part (bulk) of the spectrum of the adjacency matrix $\mathbf{A}_n$, we study its resolvent
\[
\mathbf{G}(z) \;=\; (\mathbf{A}_n/\sqrt{n} - z \mathbf{I})^{-1}, 
\]
and 
\[
S_n(z) \;=\; \frac{1}{n}\Tr\,\mathbf{G}(z) \;=\;\frac{1}{n}\sum_{i=1}^n G_{ii}(z),
\]
for $z\in\mathbb{C}^+$.  By Schur’s complement (cavity approximation), one shows
\[
G_{ii}(z)
=\Bigl(-\,z -\frac{1}{n}\sum_{j\neq i} A_{ij}\,G^{(i)}_{jj}(z)\Bigr)^{-1},
\]
where $G^{(i)}$ is the resolvent with the $i$th row and column removed.  Defining the \emph{conditional} Green’s function
\[
g_n(x;z)\;:=\;\mathbb{E}\bigl[G_{ii}(z)\mid x_i=x\bigr],
\]
and replacing each minor $G^{(i)}_{jj}(z)$ by $g_n(x_j;z)$, as well as $A_{ij}$ by its mean 
\(
\mathbb{E}[A_{ij}\mid x_i=x,\,x_j=y]
=1-e^{-n^{-1/\alpha}xy},
\)
one obtains the integral equation
\begin{equation}
g_n(x;z)=\nonumber
\Biggl(
-\,z
-
\int_{1}^{\infty}\bigl(1-e^{-n^{-1/\alpha}xy}\bigr)\,
g_n(y;z)\,\rho_\alpha(y)\,dy
\Biggr)^{-1}.
\label{eq:FPn}
\end{equation}
The above equation implicitly defines $g(x; z)$, from which the average Stieltjes transform $g_H(z) = \int \rho_\alpha(x) g(x;z),dx$ can be obtained. Solving \eqref{eq:FPn} in closed form is challenging, especially because $\rho_\alpha(x)$ has no finite moments. However, it serves as a basis for numerical solutions or asymptotic analysis. In particular, the edges of the bulk spectrum can be estimated by finding values of $z=E$ for which the denominator of \eqref{eq:FPn} approaches zero for some $x$. Intuitively, the largest bulk eigenvalue $\lambda_{\text{bulk}}^{\max}$ would satisfy $$E = \int \rho_\alpha(y) [1-e^{-\epsilon_n x E}] g(y;E)dy$$ for $x$ at the top of the weight distribution support. In practice, due to the diverging weight variance, one might expect $\lambda_{\text{bulk}}^{\max}$ to scale with $n^{1/2}$ as well (similar to the outliers).We will not delve deeper into solving \eqref{eq:FPn} here; instead, we note that a clear separation between bulk and outliers has been observed in numerical simulations (see main text). 

\subsection*{Poisson Point Process}

Since the spectrum of $\mathbf{H}$ and $\mathbf{A}$ in the bulk behave similarly (conditionally on the weights, one can show the closeness using Hoffman-Wielandt type inequality), we show that through a poisson point process approach, one can derive self-consistent equations for the matrix $\mathbf{A}$.

Let $x_{(1)}\ge x_{(2)}\ge\cdots\ge x_{(n)}$ be the order statistics of the i.i.d.\ Pareto($\alpha$) weights.  Define
\[
X_n \;=\; n^{-1/\alpha}(x_{(1)},x_{(2)},\dots,x_{(n)},0,0,\dots)\;\in\;\mathbb{R}^\infty.
\]
Classical stable-limit theorems (\cite{lepage1981convergence}) imply $X_n \xrightarrow{d}X_\infty$, where $X_\infty=(y_1,y_2,\dots)$ has
\[
y_k = \Gamma_k^{-1/\alpha}, 
\qquad
\Gamma_k = E_1 + \cdots + E_k,
\]
and $E_i\sim\mathrm{Exp}(1)$ i.i.d.  Equivalently, $\{y_k\}$ form a Poisson point process on $(0,\infty)$ with intensity
$\rho_\alpha(y)\,dy = \alpha\,y^{-1-\alpha}dy$.

Associate to each rescaled weight $w_j=x_j/n^{1/\alpha}$ the mark $h_n(w_j)=g_n(x_j;z)$.  The finite-$n$ point process
\[
N_n = \sum_{j=1}^n \delta_{(w_j,\;h_n(w_j))}
\]
converges in law to the marked Poisson process
\[
N = \sum_{k=1}^\infty \delta_{(y_k,\;g(y_k;z))},
\]
where $g(y;z)=\lim_{n\to\infty}g_n(n^{1/\alpha}y;z)$ is the \emph{limiting conditional resolvent}.
We now define the mapping
\[
\Phi_x(N_n)
\;:=\;
\sum_{j=1}^n g_n(x_j;z)\,\bigl(1 - e^{-n^{-1/\alpha}x_j\,x}\bigr),
\]
which converges to
\[
\Phi_x(N)
\;=\;
\sum_{k=1}^\infty g(y_k;z)\,\bigl(1 - e^{-x\,y_k}\bigr).
\]
The mean of $\Phi_x(N)$ is
\[
\mathbb{E}\bigl[\Phi_x(N)\bigr]
=\int_0^\infty g(y;z)\,\bigl(1-e^{-xy}\bigr)\,\rho_\alpha(y)\,dy.
\]
Hence the deterministic limit satisfies
\begin{align*}
g(x;z)
&=\Bigl(-\,z \;-\;\Phi_x(N)\Bigr)^{-1},\\
\end{align*}
or explicitly
\begin{equation}
g(x;z)
\;=\;
-\,\frac{1}
{\,z \;+\;
\displaystyle\int_0^\infty g(y;z)\,\bigl(1-e^{-xy}\bigr)\,\rho_\alpha(y)\,dy
}
\label{eq:FPinf}
\end{equation}
Finally, the bulk Stieltjes transform is obtained by averaging over the fitness distribution:
\[
g(z)
=\int_{0}^{\infty} g(x;z)\,\rho_\alpha(x)\,dx.
\]
The bulk spectral density $\rho_H(\lambda)$ follows by the Sokhotski–Plemelj formula
\[
\rho_H(\lambda)
=-\frac{1}{\pi}\lim_{\varepsilon\to0^+}\Im\,g(\lambda + i\varepsilon).
\]
Equations \eqref{eq:FPinf} and the above averaging fully characterize the continuous spectrum for $0<\alpha<1$.

\subsection*{Spectral norm bounds}

We condition on the latent variables so that the edge probabilities $(p_{ij})$ are deterministic and $\{A_{ij}\}_{i<j}$ are independent. Let
\begin{eqnarray}
v_{ij}&:=&p_{ij}(1-p_{ij}),\\ 
\sigma&:=&\max_{1\le i\le n}\Big(\sum_{j\ne i} v_{ij}\Big)^{1/2},\\
\sigma_\ast&:=&\max_{i\ne j}\sqrt{v_{ij}}\le \tfrac12.
\end{eqnarray}
Then there exist absolute constants $C,c>0$ such that
\begin{equation}\label{eq:MSM-exp}
\mathbb{E}\,\|\mathbf{H}\| \le C\Big(\sigma+\sigma_\ast\sqrt{\log n}\Big).
\end{equation}
Moreover, for any $0<\varepsilon\le \tfrac12$ and any $t\ge 0$,
\begin{equation}\label{eq:MSM-tail}
\mathbb{P}\left(\|\mathbf{H}\| \ge (1+\varepsilon)\,2\sqrt{2}\,\sigma + t\right)
\le n\exp\Big(-\,\tfrac{t^{2}}{c_\varepsilon\,\sigma_\ast^{2}}\Big),
\end{equation}
where $c_\varepsilon>0$ is universal.

\medskip

{\em Proof.}
Write $H_{ij}=\sqrt{v_{ij}}\;\xi_{ij}$ with $\xi_{ij}=\big(A_{ij}-p_{ij}\big)/\sqrt{v_{ij}}$. Conditionally on $(p_{ij})$, the $\xi_{ij}$ are independent, centered, and uniformly sub-Gaussian (bounded support), so the sub-Gaussian extension of the main bound applies with scale matrix $b_{ij}=\sqrt{v_{ij}}$. This yields \eqref{eq:MSM-exp} from the \emph{Gaussian} Theorem~1.1 together with its sub-Gaussian transfer \cite{bandeira2016sharp}{(Cor.~3.3)}: 
\begin{equation}
     \mathbb{E}\|\mathbf{H}\|\lesssim \sigma+\sigma_\ast\sqrt{\log n}.
 \end{equation}
 For tails, note that $|H_{ij}|\le 1$ and $\mathbb{E}H_{ij}=0$. Now, \cite{bandeira2016sharp}{ (Corollary~3.12)} gives, for \emph{symmetric} bounded entries, a variance-sensitive inequality of the form
\begin{equation}
    \mathbb{P}(\|\mathbf{H}\|\ge (1+\varepsilon)2\sigma+t)\le n\exp(-t^{2}/\tilde c_\varepsilon\sigma_\ast^{2}).
\end{equation}
Our entries are not symmetric, but \cite{bandeira2016sharp}{(Remark~3.13)} shows that symmetrization inflates only the leading constant by a factor $\sqrt{2}$, which yields \eqref{eq:MSM-tail}. 
As $v_{ij}\le p_{ij}$, we have $\sigma^{2}\le d_{\max}$ where $d_{\max}:=\max_i\sum_{j\ne i}p_{ij}$ is the maximal expected degree. Hence
\begin{equation}
\mathbb{E}\|\mathbf{H}\|\le C\big(\sqrt{d_{\max}}+\sqrt{\log n}\big),
\end{equation}
and
\begin{equation}
    \|\mathbf{H}\|\le C'\Big(\sqrt{d_{\max}}+\sqrt{\log n}\Big)\quad\text{w.h.p.}
\end{equation}
As a crude but universal bound, since $v_{ij}\le \tfrac14$, $\sigma\le \sqrt{n}/2$, we get

\begin{equation}
    \|\mathbf{H}\|\lesssim \tfrac{\sqrt{n}}{2}+O(\sqrt{\log n}),
\end{equation}
which can be further approximated as 
$||\mathbf{H}|| \lesssim  \frac{\sqrt{n}}{2}$
as in the main text.

\medskip

We now establish a {\em lower bound.} We condition on $(p_{ij})$, so entries are independent. Let
\begin{eqnarray}
    v_{ij} &=& p_{ij}(1-p_{ij}),\\
    \sigma^2 &=& \max_i\sum_{j\ne i} v_{ij}.
\end{eqnarray}
Then, for any $0<\delta<1$,
\begin{equation}
    \mathbb{P}\left(\,\|\mathbf{H}\|\ \ge\ \sqrt{1-\delta}\,\sigma\,\right)
\ \ge\ 1-\exp\big(-c\,\delta^2\,\sigma^2\big),
\end{equation}
for a universal constant $c>0$. In particular, if $\sigma^2\gtrsim \log n$, then $\|H\|\gtrsim \sigma$ with high probability.

\medskip

{\em Proof sketch.}
Choose $i_\star$ attaining $\sigma^2=\sum_{j\ne i_\star}v_{i_\star j}$. Since $\|H\|\ge \|He_{i_\star}\|_2$,
\begin{equation}
    \|\mathbf{H}e_{i_\star}\|_2^2=\sum_{j\ne i_\star}(A_{i_\star j}-p_{i_\star j})^2,
\end{equation}
which is a sum of independent, bounded random variables in $[0,1]$ with mean $\sigma^2$. Hoeffding’s lower-tail inequality gives
\begin{equation}
    \mathbb{P} \big(\|\mathbf{H}e_{i_\star}\|_2^2<(1-\delta)\sigma^2\big)\le \exp(-c\delta^2\sigma^2).
\end{equation}
It is predicted that $\sigma^2\gtrsim n$. Taking square roots yields the claim.

%\begin{widetext}

\section*{Eigenvectors}
\begin{figure}[h!]
    \centering 
    \includegraphics[width=0.6\linewidth]{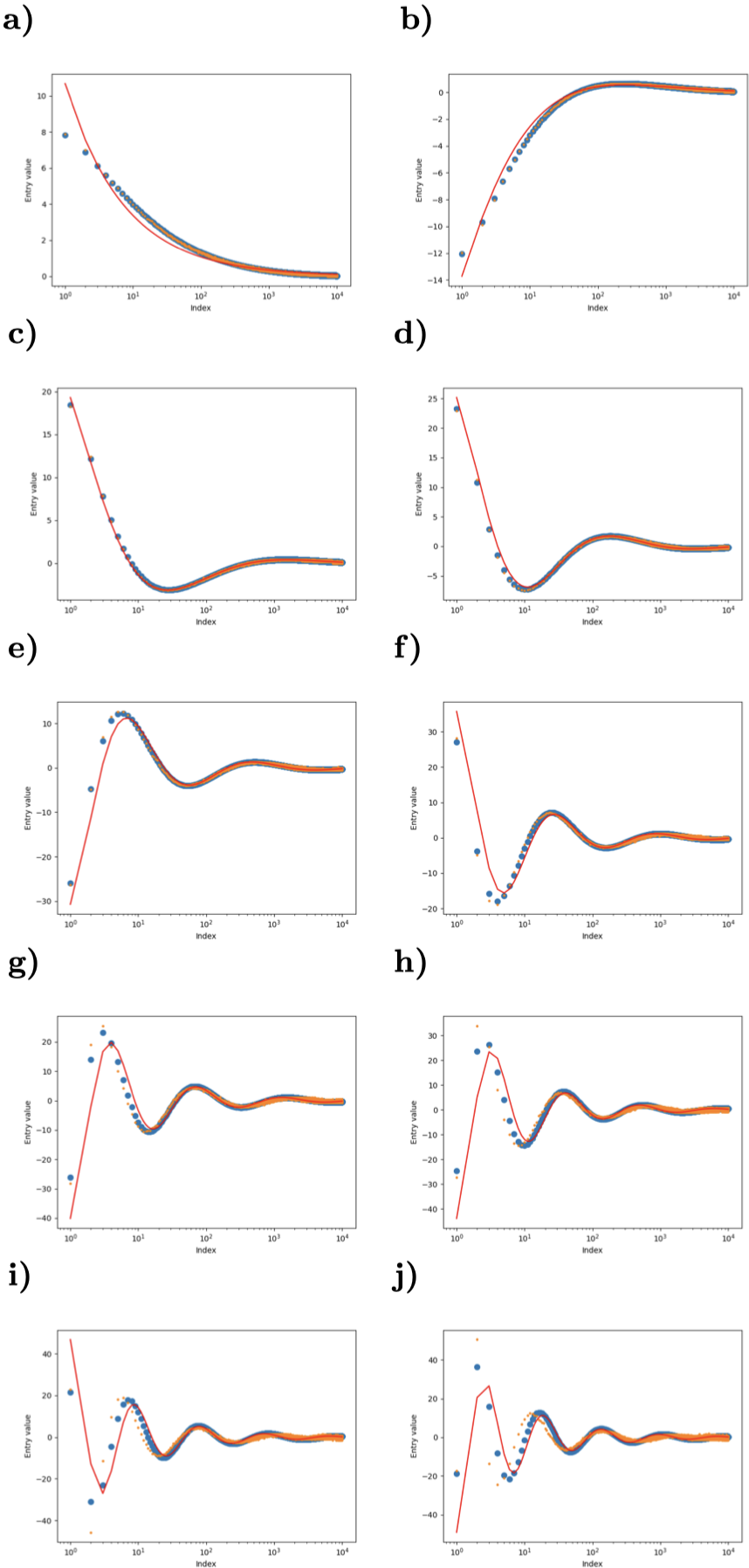}
    
    \caption{\textbf{Plot of eigenvector entries $v_k(j)$ versus $j$ for $n = 10^4$ and $\alpha = 0.2$.} \textbf{a)} First eigenvector. \textbf{b)} Second eigenvector. \textbf{c)} Third eigenvector. \textbf{d)} Fourth eigenvector. \textbf{e)} Fifth eigenvector. \textbf{f)} Sixth eigenvector. \textbf{g)} Seventh eigenvector. \textbf{h)} Eight eigenvector. \textbf{i)} Ninth eigenvector.
    \textbf{j)} tenth eigenvector. Blue dots correspond to the actual eigenvectors entries of $\mathbf{P}$, while red solid lines represent the corresponding predicted entries from \eqref{eq:eigenvec}. Orange dots are the entries of the eigenvectors of $\mathbf{A}$, and we can see that in the outlier regime they are well predicted by $\mathbf{P}$ and hence the ones we obtained analytically. }
    \label{fig:eigenvectors2}
\end{figure}

\begin{figure}[h!]
    \centering \includegraphics[width=0.6\linewidth]{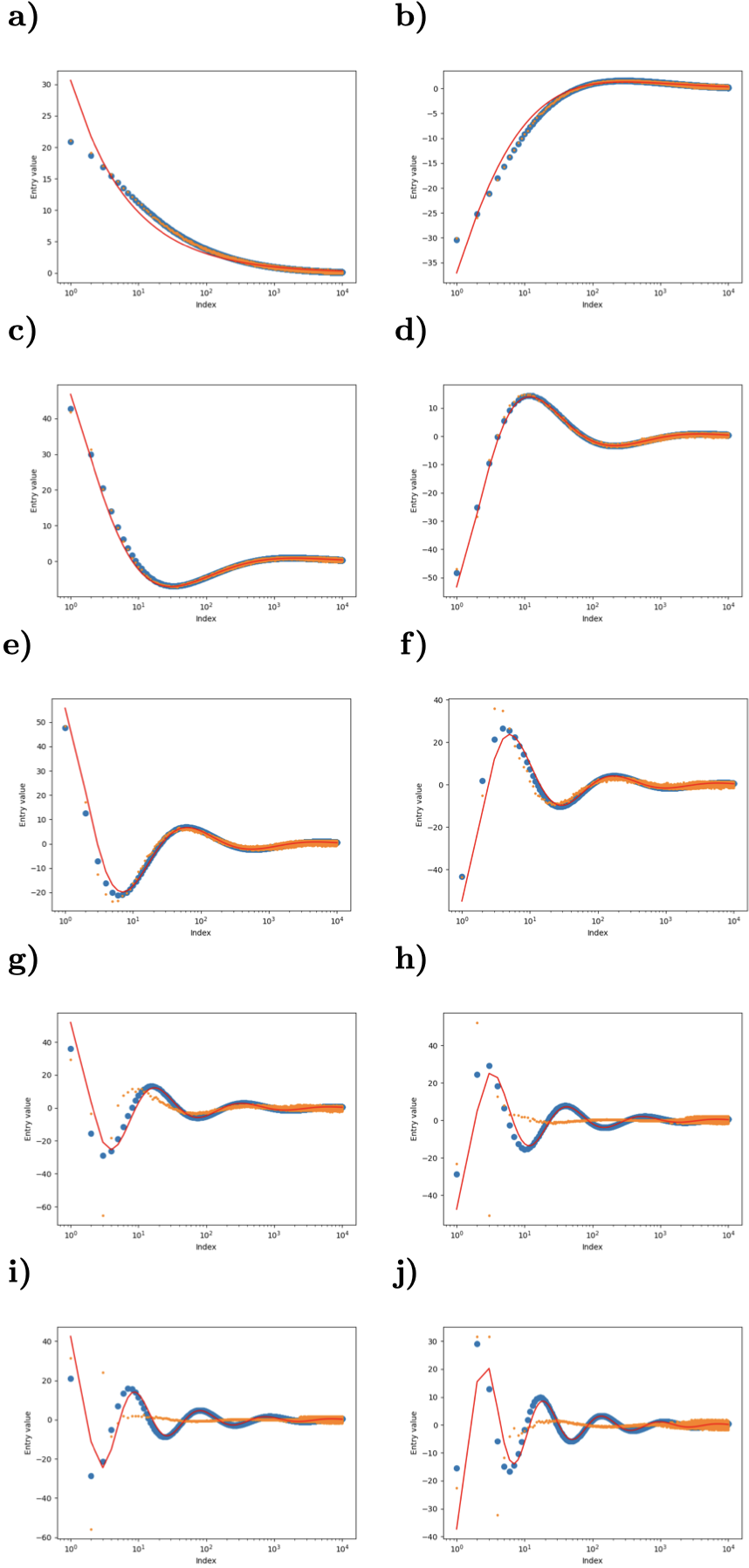}
    
    \caption{\textbf{Plot of eigenvector entries $v_k(j)$ versus $j$ for $n = 10^4$ and $\alpha = 0.5$.} \textbf{a)} First eigenvector. \textbf{b)} Second eigenvector. \textbf{c)} Third eigenvector. \textbf{d)} Fourth eigenvector. \textbf{e)} Fifth eigenvector. \textbf{f)} Sixth eigenvector. \textbf{g)} Seventh eigenvector. \textbf{h)} Eight eigenvector. \textbf{i)} Ninth eigenvector.
    \textbf{j)} tenth eigenvector. Blue dots correspond to the actual eigenvectors entries of $\mathbf{P}$, while red solid lines represent the corresponding predicted entries from \eqref{eq:eigenvec}. Orange dots are the entries of the eigenvectors of $\mathbf{A}$, and we can see that in the outlier regime they are well predicted by $\mathbf{P}$ and hence the ones we obtained analytically. }
\label{fig:eigenvectors2bis}
\end{figure}

\begin{figure}[h!]
    \centering \includegraphics[width=0.6\linewidth]{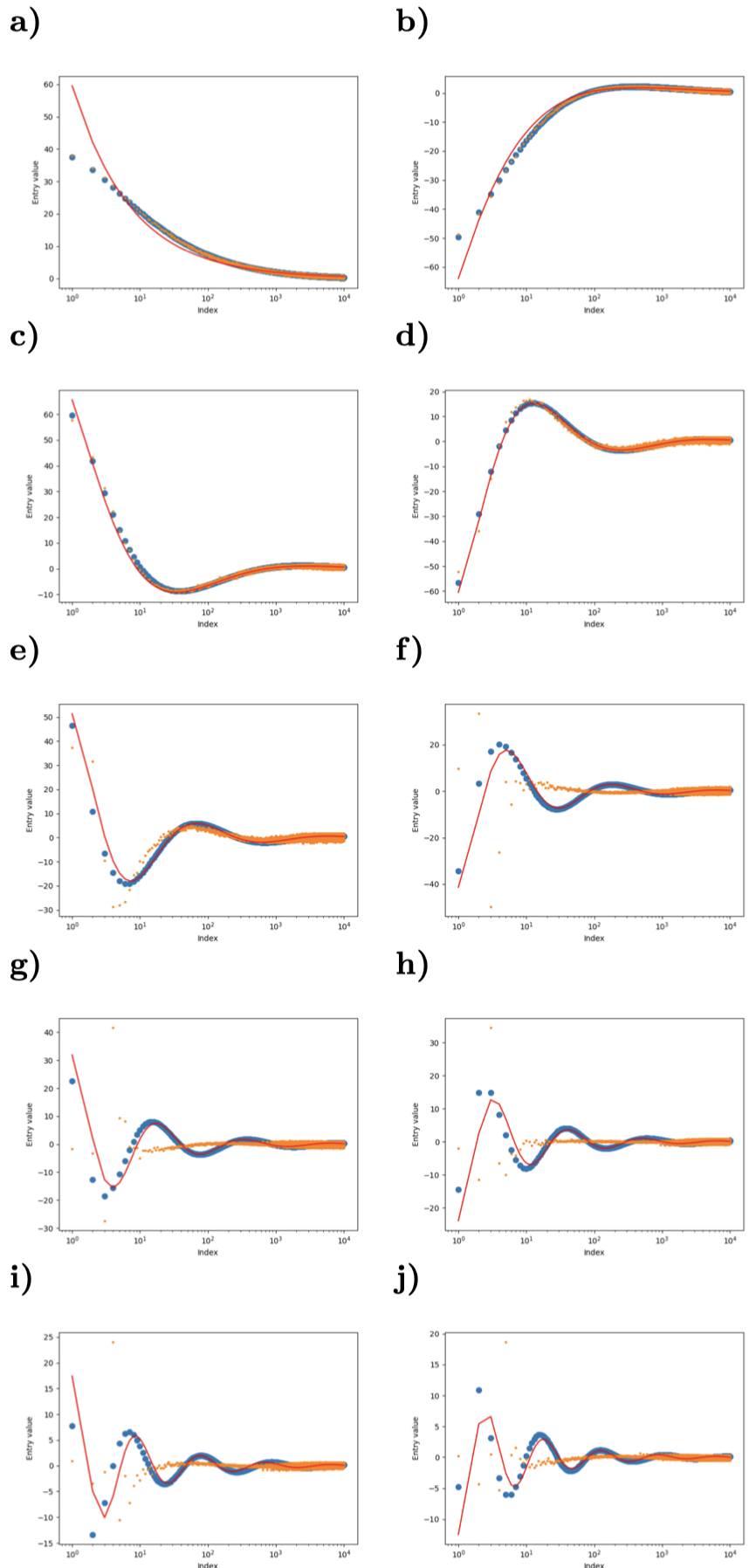}
    
    \caption{\textbf{Plot of eigenvector entries $v_k(j)$ versus $j$ for $n = 10^4$ and $\alpha = 0.8$.} \textbf{a)} First eigenvector. \textbf{b)} Second eigenvector. \textbf{c)} Third eigenvector. \textbf{d)} Fourth eigenvector. \textbf{e)} Fifth eigenvector. \textbf{f)} Sixth eigenvector. \textbf{g)} Seventh eigenvector. \textbf{h)} Eight eigenvector. \textbf{i)} Ninth eigenvector.
    \textbf{j)} tenth eigenvector. Blue dots correspond to the actual eigenvectors entries of $\mathbf{P}$, while red solid lines represent the corresponding predicted entries from \eqref{eq:eigenvec}. Orange dots are the entries of the eigenvectors of $\mathbf{A}$, and we can see that in the outlier regime they are well predicted by $\mathbf{P}$ and hence the ones we obtained analytically. }
    \label{fig:eigenvectors3}
\end{figure}

In the main text, we showed that the $j$-th entry of the eigenvector $v_k$ of the expected adjacency matrix $\mathbf{P}$ can be expressed as 
\begin{align}
    v_k^{(j)} &= \sqrt{\frac{n}{j \alpha^2}}
    ~\Re\left[\left(\frac{\alpha}{2}-\i\omega_k\right)\Gamma\left(1-\frac{\alpha}{2}-\i\omega_k\right)\left(\frac{n}{j}\right)^{\frac{\i \omega_k}{\alpha}}\right] \\
    & = -\sqrt{\frac{n}{j}}\left(\frac{\alpha}{4}+\frac{\omega_k^2}{\alpha}\right)\Re\left[\Gamma\left(-\frac{\alpha}{2}-\i\omega_k\right)\left(\frac{n}{j}\right)^{\frac{\i \omega_k}{\alpha}}\right]\label{Seq:eigenvec}
\end{align}

and that the leading eigenvectors of the realized adjacency matrix $\mathbf{A}$ follow the same structure very closely for the smallest values of $k$, but depart from it for larger values of $k$.
Here we provide additional numerical support for these results, for $\alpha=0.2$ (Fig.~\ref{fig:eigenvectors2}), $\alpha=0.5$ (Fig.~\ref{fig:eigenvectors2bis}) and $\alpha=0.8$ (Fig.~\ref{fig:eigenvectors3}).

We also note that, when combined with Eq.~\eqref{Seq:eigenvk}, the above expression allows to rewrite the eigenvector entries in terms of the corresponding eigenvalue $\lambda_k$ as
\begin{equation}
    v_k^{(j)}= \frac{(-1)^{k+1}~\lambda_k}{\sqrt{j}}   \left(\frac{\alpha}{4}+\frac{\omega_k^2}{\alpha}\right)\Re\left[~j^{\frac{\i \omega_k}{\alpha}}\right]= \frac{(-1)^{k+1}~\lambda_k}{\sqrt{j}}   \left(\frac{\alpha}{4}+\frac{\omega_k^2}{\alpha}\right)\cos{\left(\frac{\i \omega_k}{\alpha}\ln j\right)}.\label{Seq:eigenvec2}
\end{equation}
In particular, the first entry $v_k^{(1)}$ of each eigenvector $v_k$ (which is the entry that localizes around the hub node $j=1$) is
\begin{equation}
    v_k^{(1)}= {(-1)^{k+1}~\lambda_k}   \left(\frac{\alpha}{4}+\frac{\omega_k^2}{\alpha}\right),\label{Seq:eigenvec3}
\end{equation}
so that, for $j=2,n$,
\begin{equation}
    v_k^{(j)}=v_k^{(1)} \frac{\cos{\left(\frac{\i \omega_k}{\alpha}\ln j\right)}}{\sqrt{j}}.\label{Seq:eigenvec4}
\end{equation}
The above expressions illustrate how each (leading) eigenvector $v_k$ has a structure where the first entry $v_k^{(1)}$ (located at the hub node $j=1$) has an amplitude governed by the corresponding eigenvalue $\lambda_k$, and the next entries $j>1$ follow a power-law decay $\sim j^{-1/2}$ with log-periodic oscillations.

\end{document}